\documentclass[a4paper,12]{amsart}

\usepackage[left=2.5cm,top=2.5cm,right=2.5cm,bottom=2.5cm]{geometry}
\usepackage[utf8]{inputenc}
\usepackage{times} 
\usepackage{latexsym,amsfonts,amssymb,amsthm,amsmath,amscd,mathrsfs}
\usepackage{xargs}   
\usepackage{xspace}
\usepackage[shortlabels]{enumitem}
\usepackage{graphicx}
\usepackage{float}
\usepackage{multirow}
\usepackage[pdftex,dvipsnames]{xcolor} 
\usepackage{hyperref}
\usepackage{graphicx,tikz,times,pgffor}
\usepackage{pgfplots}
\usetikzlibrary{calc,through,backgrounds,patterns}
\usepackage{tabulary}
\usepackage{tabularx}
\usepackage{soul}
\usepackage{caption}
\usepackage{subcaption}
\usepackage{mathtools}
\usepackage{comment}

\pgfplotsset{compat=newest} 
\pgfplotsset{plot coordinates/math parser=false} 
\newlength\figureheight 
\newlength\figurewidth 

\usepackage{xpatch}

\usepackage{algpseudocode}
\usepackage{algorithm}

\newcommand{\testColor}{Black}
\newcommand{\setB}{\mathcal{B}}
\newcommand{\setC}{\mathcal{C}}

\newcommand{\setFT}{{\color{\testColor}\setT_{\mathcal{A}}}}
\newcommand{\setFX}{{\color{\testColor}\setX_{\mathcal{A}}}}
\newcommand{\setFXapprox}{\color{\testColor}{{\color{\testColor}\widehat\setFX}}}
\newcommand{\setP}{\mathcal{P}}
\newcommand{\setT}{{\color{\testColor}\mathbb{T}}}
\newcommand{\setTs}{{\color{\testColor}\setT_s}}
\newcommand{\setTinv}{{\color{\testColor}\mathbb{T}_\mathrm{inv}}}
\newcommand{\setTinvD}{{\color{\testColor}\mathbb{T}^d_\mathrm{inv}}}
\newcommand{\setTA}{{\color{\testColor}\setT_{\mathcal{A}}}}
\newcommand{\setTAapprox}{{\color{\testColor}\widehat\setT_{\mathcal{A}}}}
\newcommand{\setTAinv}{\ensuremath{{\color{\testColor}\mathbb{T}_\mathrm{inv}}}}
\newcommand{\setU}{{\color{\testColor}\mathbb{U}}}
\newcommand{\setX}{{\color{\testColor}\mathbb{X}}}
\newcommand{\setXs}{{\color{\testColor}\setX_s}}
\newcommand{\setXsD}{{\color{\testColor}\setX_s^d}}
\newcommand{\setXinvD}{{\color{\testColor}\mathbb{X}^d_\mathrm{inv}}}
\newcommand{\setXinv}{{\color{\testColor}\mathbb{X}_\mathrm{inv}}}
\newcommand{\setXA}{{\color{\testColor}\setX_{\mathcal{A}}}}
\newcommand{\setY}{{\color{\testColor}\mathbb{Y}}}
\newcommand{\setYA}{{\color{\testColor}\setY_{\mathcal{A}}}}
\newcommand{\setYinv}{{\color{\testColor}\mathbb{Y}_\mathrm{inv}}}
\newcommand{\setYs}{{\color{\testColor}\mathbb{Y}_{s}}}
\newcommand{\setYsD}{{\color{\testColor}\mathbb{Y}_{s}^d}}

\newcommand{\setZ}{{\color{\testColor}\mathcal{Z}}}
\newcommand{\Ucal}{{\color{\testColor}\mathcal{U}}}

\newcommand{\setPolyXinv}{\ensuremath{{\color{\testColor}\setFXapprox^\mathrm{inv}}}}
\newcommand{\setPolyXinvT}{\ensuremath{{\color{\testColor}\setTAapprox^\mathrm{inv}}}}

\newcommand{\dist}{\mathrm{dist}}
\DeclareMathOperator*{\argmin}{arg\,min}

\newcommand{\R}{\mathbb{R}}
\newcommand{\N}{\mathbb{N}}

\DeclareMathOperator{\inti}{int}

\newtheorem {assumpt} {Assumption}
\newtheorem {theorem}{Theorem}
\newtheorem {defi}{Definition}
\newtheorem{rem}{Remark}
\newtheorem{lem}{Lemma}
\newtheorem {propt}{Property}

\newtheorem {cor}{Corollary}
\newtheorem* {probstat} {Control Problem}

\newenvironment{pf*}[1][Proof:]{\noindent \textbf{#1} }{}
\newenvironment{spf}{\noindent \textbf{Sketch of Proof:}}{}

\usepackage[shadow,colorinlistoftodos,prependcaption,spanish,textsize=footnotesize]{todonotes}

\newcommandx{\falta}[2][1=]{\todo[linecolor=red,backgroundcolor=red!25,bordercolor=red,#1]{#2}}
\newcommandx{\completar}[2][1=]{\todo[linecolor=blue,backgroundcolor=blue!25,bordercolor=blue,#1]{#2}}
\newcommandx{\chequear}[2][1=]{\todo[linecolor=OliveGreen,backgroundcolor=OliveGreen!25,bordercolor=OliveGreen,#1]{#2}}
\newcommandx{\improve}[2][1=]{\todo[linecolor=Plum,backgroundcolor=Plum!25,bordercolor=Plum,#1]{#2}}
\newcommandx{\explain}[2][1=]{\todo[linecolor=lime,backgroundcolor=lime!25,bordercolor=lime,#1]{#2}}
\newcommandx{\marce}[2][1=]{\todo[linecolor=magenta,backgroundcolor=magenta!25,bordercolor=magenta,#1]{#2}}
\newcommandx{\igna}[2][1=]{\todo[linecolor=SeaGreen,backgroundcolor=SeaGreen!25,bordercolor=magenta,#1]{#2}}
\newcommandx{\ig}{\textcolor{SeaGreen}}

\numberwithin{equation}{section}

\allowdisplaybreaks
\title[Characterization and computation of control invariant sets within target regions for linear ICS]{Characterization and computation of control invariant sets within target regions for linear impulsive control systems}
\author[I. Sanchez, C. Louembet, M. Actis and A. H. Gonzalez]{Ignacio Sanchez$^\dagger$, Christophe Louembet$^{\ddagger}$, Marcelo Actis$^{\ast}$ and Alejandro H. Gonzalez$^{\dagger\dagger}$}
\thanks{$^{\dagger}$Institute of Applied Mathematics of Litoral (IMAL), CONICET-UNL, Santa Fe, Argentina (e-mail: isanchez@santafe-conicet.gov.ar)}
\thanks{$^{\ddagger}$LAAS-CNRS, Universit\'e de Toulouse, CNRS, Toulouse, France (e-mail: louembet@laas.fr)}
\thanks{$^{\ast}$Facultad de Ingeniera Qumica (FIQ), Universidad Nacional del Litoral (UNL) and Consejo Nacional de Investigaciones cientficas y tecnicas (CONICET), Santa Fe, Argentina}
\thanks{$^{\dagger\dagger}$Institute of Technological Development for the Chemical Industry (INTEC), CONICET-UNL, Santa Fe, Argentina (e-mail: alejgon@santafe-conicet.gov.ar)}
\date{\textcolor{black}{\today}}

\begin{document}

\maketitle

\begin{abstract}
Linear impulsively controlled systems are suitable to describe a venue of real-life problems, going from disease treatment to aerospace guidance. 
The main characteristic of such systems is that they remain uncontrolled for certain periods of time.
As a consequence, punctual equilibria characterizations outside the origin are no longer useful, and the whole concept of equilibrium and its natural extension, the controlled invariant sets, needs to be redefined. 
Also, an exact characterization of the admissible states, 
i.e., states such that their uncontrolled evolution between impulse times remain within a predefined set, is required.
An approach to such tasks --- based on the Markov-Lukasz theorem --- is presented, providing a tractable and non-conservative  characterization, emerging from polynomial positivity that has application to systems with rational eigenvalues. This is in turn the basis for obtaining a tractable approximation to the maximal admissible invariant sets. 
In this work, it is also demonstrated that, in order for the problem to have a solution, an invariant set (and moreover, an equilibrium set) must be contained within the target zone.
To assess the proposal, the so-obtained impulsive invariant set is explicitly used in the formulation of a set-based model predictive controller, with application to zone tracking.
In this context, specific MPC theory needs to be considered, as the target is not necessarily stable in the sense of Lyapunov. A zone MPC formulation is proposed, which is able to i) track an invariant set such that the uncontrolled propagation fulfills the zone constraint at all times and ii) converge asymptotically to the set of periodic orbits completely contained within the target zone.
\end{abstract}


\section{Introduction}

Impulsive systems are a subclass of dynamical systems (indeed, a hybrid dynamical system) in which resetting or impulsive events produce a discontinuity of the first kind in the state trajectories. This kind of systems has been extensively studied in the literature, and results concerning the existence, uniqueness and stability of solutions have been achieved (\cite{bainov1995,Yang01,haddad2006}). 

Particularly, when the impulsive nature of the control actions (inputs) is the one that produces the discontinuity, we have an impulsively controlled system. Many control problems fall in the scope of impulsively controlled systems, as it is the case of drug scheduling in several disease treatments (by taking pills \cite{Rivade2019,rivadeneira2018non,Riv2017,magdaleno2019learning,hernandez2019passivity} or by applying injections \cite{AbuinJPC20,GonOCAM20}), or that of aircraft guidance \cite{ArantesPHD2018,Arantes2019,louembetCDC2019impulsive}.
The objective in this kind of control problem is to maintain the closed-loop system in a target region (defined by the problem itself) where the operation is safe. According to the meaning of the state variables, these regions does not include the origin as an interior point, since the origin uses to be a state of emptiness or rest, which is not safe. So, there is no formal equilibrium in the target region, and the control objectives are redefined to steer the closed-loop system as close as possible to the target region. 

To accomplish the control objective properly ---i.e., to ensure that the state will reach and remain inside the target region --- it is necessary to formally define both, an extended equilibrium and extended controlled invariant set.
Following the ideas in \cite{SopaPatriBempo15} and \cite{RivTAC18} a two-set definition can be used, in which an equilibrium or controlled invariant set is such, only with respect to a larger set that contains the state free responses between the impulses. 
Although potentially conservative (since the outer set can be as large as desired) these definitions show to be useful to define equilibrium and controlled invariant sets with respect to target regions, and provides formal tools to set the control problem, mainly when the controller is a model-based one \cite{SopaPatriBempo15,RivTAC18}. 
In addition, based on this two-set generalization, formal notions of impulsive closed-loop convergence and stability \cite{Lobo_Pereira_2015,djorge2020stability} can be applied.

In this work, we discuss the way to formally compute exact admissible sets for impulsively controlled linear systems, with respect to a given target set. 
In the seminal work \cite{gilbert1991linear}, the study of admissible sets for linear continuous-time system is approached, where the initial state of an unforced linear system is called output admissible with respect to a constraint set $\setY$ if the resulting output function satisfies the pointwise-in-time condition $y(t) \in \setY, t \geq 0$. In this work, we will refer to admissible sets for constrained state trajectories during the interval between impulses. For sake of shortness, we will omit to denote with admissible certain sets as equilibrium, feasible and invariant sets, where their admissibility will be evident from context.
This notion, already introduced into the analysis of linear impulsive systems in \cite{SopaPatriBempo15} and mentioned in \cite{RivTAC18}, was not exactly characterized yet: only approximating techniques are proposed. 
In such context, a first contribution of this paper is to solve - by means of a tractable methodology - the particular problem of characterizing the set of states for which the corresponding free trajectories remain in a given target set (this set denoted as admissible set). 
The free propagation of the states is expressed as univariate polynomial thanks to a relevant change of variable, and then the set of trajectories included in the target set is characterized through the Linear Matrix Inequalities (LMI) conditions on the initial states.
These conditions are based on the Markov-Lukasz theorem (see \cite{Nesterov2000}) and they have been inspired by the work \cite{Henrion_2005}. 

A procedure is then proposed to compute the maximal invariant set with respect to the admissible set. The resulting set - defined as impulsive invariant set - satisfies the constraints at all times, meaning property that the system state remains within the set at the impulsive times (the times at which the impulsive inputs enter the system) and does not leave the target set at any time. 

A novel set-based MPC algorithm is proposed, that exploits the tractability and utility of this description, and follows the procedure exposed in \cite{AndersonOCAM18,GonzalezSCL2014,FerramoscaJPC10}
extended through the explicit use of the new characterizations of invariant sets. This way, a closed-loop stable controller able to steer the system to the target region, and to maintain it in such a set indefinitely was obtained, with guaranteed feasibility at all times.

The main theoretical contribution of this work is the development of a proof of necessary conditions for the validity of a target zone for a zone tracking control with impulsive control inputs. 

The effectiveness of the proposed controller formulation is assessed through simulation.

\subsection{Notation}\label{sec:notation}

Let $\setX \subseteq\R^n$ and $\setY\subseteq\R^m$. A \emph{correspondence} $c: \setX \rightrightarrows \setY$ defines for each $x \in \setX$ a set
$c(x) \subseteq \setY$. A correspondence (also denoted as \emph{set-valued function}) is a generalization of the concept of \emph{function}, $f$, which
for each $x \in \setX$ defines a unique $f(x) \in \setY$.
The \emph{euclidean distance} between two points $x,y \in\R^n$ is denoted by $\|x-y\|:=[(x-y)'(x-y)]^{1/2}$. 
The \emph{distance} from $x \in \R^n$ to $\setX \subseteq \R^n$ is given by $\mbox{Dist}_{\setX}(x):= \mbox{inf}_{y\in\setX} \|x-y\|$. 
The Minkowski sum $\setX \oplus \setY$ is defined by $\setX \oplus \setY := \{x+y : x\in\setX,y\in\setY\}$.
The \emph{open ball with center in $x\in \setX$ and radius $\varepsilon>0$} is given by {$\mathcal{B}_{\varepsilon}(x):=\{y \in \setX : \|x- y\|<\varepsilon\}$}.
The \emph{$\varepsilon$-neighborhood of set $\setX$} is given by $\setB_{\varepsilon}(\setX):=\{\setX \oplus \setB_{\varepsilon}(0)\}$.
Given $x\in \setX$, we say that $x$ is an \emph{interior point of $\setX$} if it there exists
$\varepsilon>0$ such that the open ball $\mathcal{B}_{\varepsilon}(x) \subseteq \setX$. The \emph{interior of $\setX$} is the 
set of all interior points and it is denoted by $\inti \setX$.
The infinite sequence of elements, $\{x_k\}_{k=1}^{\infty}$, will be simply denoted as $\{x_k\}$.

\section{Preliminaries}

First, consider the following impulsively controlled linear system (ICSys)
\begin{eqnarray}\label{ec:sysIm}
\dot{x}(t)  & = & A x(t),~~~~~~~~~~~~~~~~~~~~~~~~~~~ t \neq \tau_k, \\ 
x(\tau_k) & = & x(\tau^-_k) +B u(\tau_{k-1}),~~~~~ k \in \mathbb{N}, \nonumber
\end{eqnarray}
where $x \in \setX \subset \R^n$ represents the state, $A \in \R^{n\times n}$ is the transition matrix, $u \in \setU 
\subset \R^m$ is the input, $\tau_k =kT$, for a time period $T>0$ is the jump or impulse time, and $\tau_k^-$ denotes the time just before $\tau_k$ (i.e.,
$x(\tau_k^-) = \lim_{\delta \rightarrow 0^+} x(\tau_k - \delta))$.
The state set $\setX$ is assumed to be a closed polyhedron, the input set $\setU$ is assumed to be a compact polyhedron, and both are assumed to 
contain the origin in their nonempty interior.

For any $t \in [\tau_k,\tau_{k+1})$, $ k \in \mathbb{N}$, the solution of \ref{ec:sysIm} is given by $x(t)=e^{At}x(\tau_k)$, and $x(\tau_{k+1})=x(\tau^-_{k+1}) + B u(\tau_k)$. Given that $x(\tau^-_{k+1})=e^{AT}x(\tau_k)$, then, we can write
\begin{eqnarray}\label{ec:sol}
x(t)  & = & e^{At}x(\tau_k),~~~~~ t \in [\tau_k,\tau_{k+1}),\\
x(\tau_{k+1})  & = & e^{AT}x(\tau_k) + B u(\tau_k), \nonumber
\end{eqnarray}
which describe the free response and the jump produced by the input, respectively.

In order to properly characterize the equilibrium and invariant sets in the next sections, the following definitions concerning the uncontrolled or free responses is made.
\begin{defi}[Admissible set]\label{defi:Tset}
	Consider the ICSys system \eqref{ec:sysIm} and a polytopic non-empty set $\setY \subseteq \setX$. The admissible set of $\setY$ is given by
	\begin{eqnarray*}
		\setYA := \{x \in \setY:~e^{At}x \in \setY,~t\in [0,T] \}.
	\end{eqnarray*}
\end{defi}
Set $\setYA$ is the set of initial states for which the free responses - that are independent of $u$ - remain in $\setY$ for at least an interval of length $T$. As in \cite{SopaPatriBempo15}, this set can be described by the intersection of an uncountable set of constraints as follows:
\begin{equation}
    \setYA := \bigcap_{\tau \in [0,T)} e^{-A \tau} \setY.
\end{equation}
Note that for each $\tau \in [0,T)$ the set $e^{-A \tau} \setY$ is polytopic, since $e^{-A \tau}$ is a linear map, so $\setYA$ is given by the intersection of (an uncountable number of) polytopes. As a result, set $\setYA$ is closed and convex and can be characterized and computed by following the ideas presented in \cite{Nesterov2000,henrion2005control,Arantes2018} concerning spectrahedron representation. The spectrahedron representation of admissible sets is one of the main contribution of this work, that allows us to compute equilibrium and invariant sets for system \eqref{ec:sysIm}, as detailed next in Subsection \ref{sec:AdmCharact}.

\subsection{Sampling the ICS}
%
By sampling the ICSys \eqref{ec:sysIm} at times $\tau_k$, $k \in \mathbb{N}$, the Discretized ICSsys (DICSys) ---i.e., a discrete-time system associated to the ICSys--- is obtained \footnote{Note that this is the sampling for an impulsively controlled system, which is different from sampled-data formulations, where usually a zero order hold (ZOH) is assumed for the inputs.}
\begin{eqnarray} \label{eq:dsys}
x(\tau_{k+1}) = A^d x(\tau_{k}) + B^d u(\tau_{k}), \label{ec:sysDisc}
\end{eqnarray}
with $A^d = e^{AT}$ and $B^d = B$. The idea is to use this simplified system to infer properties of
the ICSys. Particularly, easy computations of controlled invariant sets for the ICSys will be obtained
based on controlled invariant sets for DICSys.
To this end, and for the sake of clarity, the following definitions ---concerning classical equilibria and controlled invariant sets
for discrete-time systems--- are recalled. 
\begin{assumpt}
The pair $(A^d,B^d)$ is controllable and the state is measured at each sampling time.
\end{assumpt}
\begin{defi}[Controlled equilibrium set (CES)]\label{defi:ces}
	Consider the DICSys system \eqref{ec:sysDisc}. A nonempty convex set $\setXsD \subset \setX$ is a controlled equilibrium set if for every $x_s \in \setXsD$ exists $u_s \in \setU$ such that $x_s = A^d x_s + B^d u_s$.
\end{defi}
\begin{defi}[Controlled invariant set (CIS) \cite{Blanchinibook15}]\label{defi:cis}
	Consider the DICSys system \eqref{ec:sysDisc}. A nonempty convex set $\setXinvD \subset \setX$ is a controlled invariant set	if for every $x \in \setXinvD$ exists $u \in \setU$ such that $A^d x + B^d u \in \setXinvD$. A CIS with nonempty interior is denoted a proper CIS.
\end{defi}
Clearly, every CES is a CIS (although exceptionally are a proper CIS, since they use to have an empty interior). However, the fact that any CIS contains a CES is not trivial. In \cite{feuer1976omega} (Theorem 3.3) a proof is given for the continuous-time case, based on the theorem of Kakutani \cite{kakutani1941generalization} (a generalization of the Brouwer's fixed-point theorem). The next theorem provides a similar result for linear discrete-time systems with polytopic constraints, directly based on the Brouwer's fixed-point theorem (Theorem~\ref{theo:fixPoint} in Appendix \ref{sec:AppA}).
\begin{theorem}\label{theo:Main}
	Consider the DICSys \eqref{ec:sysDisc}. Then, every compact and convex CIS $\setXinvD \subset \setX$ contains a CES $\setXsD$ (this set may be a singleton $\setXsD = \{x_s\}$).
\end{theorem}

\begin{proof}
We present here a new proof of Theorem \ref{theo:Main} (different from the one given in \cite{feuer1976omega}) that is directly based on the Brower fixed point.
Consider a compact and convex CIS $\setXinvD$ and the autonomous (or closed-loop) system $x^+=f_{\kappa}(x)$, where $f_{\kappa}(x):= A^dx + B^d\kappa(x)$, for each
	$x\in \setXinv$, and $\kappa(x)$ is defined by
	\begin{eqnarray}\label{ec:u_of_x}
	\kappa(x) := \argmin_{u \in \Ucal(x)} V(x,u),
	\end{eqnarray}
	being $\Ucal(x) := \{ u \in \setU : A^dx + B^du \in \setXinvD\}$ the set of all $u$'s that keep a particular $x \in \setXinvD$
	in $\setXinvD$, and $V(x,u)$ a real convex function on $\setU$, for each $x\in \setXinvD$ (for instance, $V(x,u):=u^2$).
		
	In the optimization problem \eqref{ec:u_of_x}, $u$ is the optimization/decision variable, $x$ the optimization parameter
	and $\Ucal(x)$ is a correspondence ($\Ucal: \setXinvD \rightrightarrows \setU$). By the definition of CIS,
	the correspondence $\Ucal(x)$ \textbf{is non-empty for each} $x \in \setXinvD$, and $\Ucal(x)$
	\textbf{is convex} because $\setXinvD$ is convex and DICSys \eqref{ec:sysDisc} is linear. Furthermore,
	$\Ucal(x)$ \textbf{is compact} because both, $\setXinvD$ and $\setU$ are compact and the DICSys \eqref{ec:sysDisc} is linear.
	
	To show that $\setU(x)$ \textbf{is continuous} on $\setXinvD$, it is necessary to show that it is both upper and lower semicontinuous (according to Definitions \ref{defi:upsemi} and \ref{defi:lowsemi}, in Appendix \ref{sec:AppB}), which is proved in Lemmas \ref{lem:uppercontin} and \ref{lem:lowercontin}, respectively, in the same Appendix \ref{sec:AppB}.
	
	Summarizing, we have that $V(x,u)$ is convex (and so, strictly cuasi-convex) and continuous on $\setU \supset \Ucal(x)$, while $\Ucal(x)$ is non-empty, convex, compact for each $x \in \setXinvD$, and continuous on $\setXinvD$. Then, by Theorem \ref{theo:VarBerge}, in Appendix \ref{sec:AppB}, it follows that $\kappa(x)$ is a continuous function on $\setXinvD$, which implies that $f_{\kappa}(x)$ is also a continuous function.
	Finally, it is easy to see that $f_{\kappa}$ maps the compact and convex CIS $\setXinvD$ into itself and so, by the Brouwer's fixed point theorem (Theorem \ref{theo:fixPoint}, in Appendix \ref{sec:AppA}), there is a state $x_s \in \setXinvD$ such that $f_{\kappa}(x_s)=x_s$. This means that $x_s =A^dx_s + B^d\kappa(x_s)$, i.e., for the fixed point $x_s$ there exists an $u_s \in \setU$, $u_s:=\kappa(x_s)$, such that 
	$x_s$ is a controlled equilibrium point for DICSys \eqref{ec:sysDisc}, which conclude the proof.
\end{proof}

\begin{rem}
	Note that the compactness and convexity of $\setXinvD$ are crucial in Theorem \ref{theo:Main}. It is quite easy to find
	counterexamples otherwise (see \cite{Blanchinibook15}, chapter 4, exercise 3).
\end{rem}

\subsection{Controlled equilibrium and invariant sets for ICSys}
%
In the context of ICSys, the equilibrium and its generalizations needs to be defined in a quite general form in contrast to typical continuous or discrete-time systems. 
Indeed, given that there are periods of uncontrolled state evolution, the following definitions are necessary:
\begin{defi}[Impulsively controlled equilibrium set (ICES) \cite{RivTAC18}]\label{defi:equilset}
	Consider the ICSys system \eqref{ec:sysIm} and a convex set $\setY \subseteq \setX$. A nonempty convex set $\setYs \subseteq \setY$ is an impulsive controlled equilibrium set if for every $x_s \in \setYs$, it follows that: (i) $\{e^{At}x_s:~t \in [0,T]\} \subset \setY$ and (ii) it there exists $u_s \in \setU$ such that $e^{AT}x_s+Bu_s = x_s$. Every single state $x_s$ in $\setYs$ is denoted as impulsive equilibrium state w.r.t $\setY$.
\end{defi}

\begin{rem}
	Note that in general, the only formal controlled equilibrium pair $(x_s,u_s)$ of ICSys \eqref{ec:sysIm} --- i.e., those that requires that both the jumps and the free response remains in a fixed state --- is the origin, i.e.,$(x_s,u_s)=(0,0)$ \cite{Yang01}, since this is the unique no-jump scenario. This fact critically shortens the application scope of impulsive control systems representations, since the origin (meaning rest or emptiness) is not included into the target sets (zones) of interest \cite{SopaPatriBempo15,RivTAC18}.
\end{rem}
\begin{defi}[Impulsive controlled invariant set (ICIS) \cite{SopaPatriBempo15}]\label{defi:invset}
	Consider the ICSys system \eqref{ec:sysIm} and a convex set $\setY \subseteq \setX$. A nonempty convex set $\setYinv \subset \setY $ is an impulsive controlled invariant set if for every $x \in \setYinv$ exists $u \in \setU$ such that (i) $\{e^{At}x:~t \in [0,T]\} \subset \setY$ and (ii) $e^{AT}x+Bu \in \setYinv$.	An ICIS with nonempty interior is denoted a proper ICIS.
\end{defi}

Note that neither $\setYs$ nor $\setYinv$ are unique for a given $\setY$. Clearly, for a given $\setY$, an ICES is an ICIS (although usually it is not a proper ICIS, i.e, it has an empty interior). 
Furthermore, ICES and ICIS for ICISys \eqref{ec:sysIm} are also CES and CIS for the DICSys \eqref{ec:sysDisc} 
(i.e., $e^{AT}(x+Bu)=A^d x+B^d u$), and CES and CIS for the DICSys \eqref{ec:sysDisc} are ICES and ICIS for ICISys \eqref{ec:sysIm} if $\setY = \R^n$. 
Finally, note also that an ICIS, $\setYinv$, can be seen as a CIS contained in $\setYA$ (this subtle fact will be used next to characterize and compute ICIS).

A question that naturally arises at this point is if Theorem \ref{theo:Main} can be extended to the case of ICISys \eqref{ec:sysIm}.
\begin{theorem}\label{theo:MainImp}
	Consider the ICSys \eqref{ec:sysIm} and a convex set $\setY \subseteq \setX$. Then, every compact and convex ICIS, $\setYinv$, contains an ICES, $\setYs$.
\end{theorem}
\begin{proof}
	Note first that conditions (ii) in Definitions \ref{defi:equilset} and \ref{defi:invset} mean that the ICES and ICIS w.r.t. $\setY$, for the ICSys \eqref{ec:sysIm}, are also CES and CIS, respectively, for the corresponding DICSys \eqref{ec:sysDisc}. Then, $\setYinv$ is also a CIS for the DICSys and contains a CES, by Theorem \ref{theo:Main}. Denote this CES as $\setYsD$. Given that any state $x_s$ in $\setYsD$ is also in $\setYinv$, and $\setYinv$ is
	an ICIS w.r.t. $\setY$, then $\{e^{At}x_s,~t\in [0,~T]\} \in \setY$, which means that $\setYs$ is a ICES w.r.t. $\setY$.
\end{proof}
\subsection{Control Problem}\label{sec:contrprob}
A control problem that frequently arises in applications of ICSys, known as Zone Tracking,  can be stated as follows:
\begin{probstat}
	Given a compact convex set $\setT \subseteq \setX$, denoted as the \textbf{target set}, coming from the application itself, the control objective is to feasibly steer the system state to $\setT$, and maintain it there indefinitely. 
\end{probstat}
To feasibly steer the system to somewhere means fulfill the input and state constraints at all time instants, so the control problem is closed related to the admissible set of $\setT$ and $\setX$, $\setTA$ and $\setXA$, respectively. Furthermore, to remain indefinitely in a given set, by means of feasible control actions, means invariance.

Therefore, according to this latter concept, the set $\setT$ coming directly from the applications needs to be refined as
%

\begin{defi}[Valid target set]\label{defi:targ}
	Consider the ICSys system \eqref{ec:sysIm}. A target set $\setT \in \setX$ (coming from the control problem definition) is a 
	\textbf{valid target set} if it contains an impulsive controlled invariant set (ICIS)(i.e., it there exists a non-empty $\setTinv$).
\end{defi}

The latter definition makes sense for most of the applications. Indeed, if the conditions of Definition \ref{defi:targ} are not fulfilled, then the problem is not well formulated, since it is not possible to remain in $\setT$ indefinitely.
Note that, in some cases, to ensure that a given $\setT$ is a valid objective region it is necessary to reduce the time period $T$ since shorter time periods produce smaller \emph{drifts} (see \cite{RivTAC18} for details).

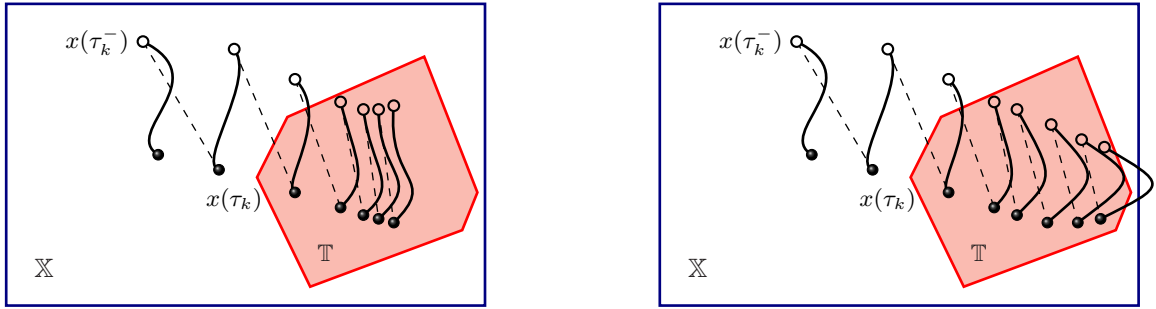
\begin{figure}[!t]
\begin{subfigure}[b]{0.45\textwidth}
	\centerline{
	\begin{tikzpicture}[scale=1]
	\draw[blue!50!black, line width=1.0pt] (3,3.5) --  (9.3,3.5) -- (9.3,-0.5) -- (3,-0.5) -- cycle;
	\draw[red, line width=1.0pt, fill=Salmon!60] (7,-0.25) --  (6.3,1.2) -- (6.7,2) -- (8.5,2.8) -- (9.2,1) -- (9,0.5) -- cycle;
	\draw[black, line width=1.0pt] (5,1.5) .. controls (4.5,1.8) and (5.8,2.5) .. (4.8,3);
	\draw[black, line width=0.5pt, dashed] (4.8,3) -- (5.8,1.3); 
	\draw[black, line width=1.0pt] (5.8,1.3) .. controls (5.5,1.4) and (6.3,2.5) .. (6,2.9);
	\draw[black, line width=0.5pt, dashed] (6,2.9) -- (6.8,1); 
	\draw[black, line width=1.0pt] (6.8,1) .. controls (6.6,1.1) and (7.3,2) .. (6.8,2.5);
	\draw[black, line width=0.5pt, dashed] (6.8,2.5) -- (7.4,0.8); 
	\draw[black, line width=1.0pt] (7.4,0.8) .. controls (7.8,1.2) and (7.6,1.4) .. (7.4,2.2);
	\draw[black, line width=0.5pt, dashed] (7.4,2.2) -- (7.7,0.7); 
	\draw[black, line width=1.0pt] (7.7,0.7) .. controls (8.2,1.1) and (7.8,1.3) .. (7.7,2.1);
	\draw[black, line width=0.5pt, dashed] (7.7,2.1) -- (7.9,0.65); 
	\draw[black, line width=1.0pt] (7.9,0.65) .. controls (8.5,1.2) and (7.9,1.2) .. (7.9,2.1);
	\draw[black, line width=0.5pt, dashed] (7.9,2.1) -- (8.1,0.6); 
	\draw[black, line width=1.0pt] (8.1,0.6) .. controls (8.7,1.2) and (8,1.1) .. (8.1,2.15);
	\draw[ball color=black] (5,1.5) circle (2 pt); 
	\draw[black, thick,fill=white] (4.8,3) circle (2 pt);
	\draw[ball color=black] (5.8,1.3) circle (2 pt); 
	\draw[black, thick,fill=white] (6,2.9) circle (2 pt);
	\draw[ball color=black] (6.8,1) circle (2 pt); 
	\draw[black, thick,fill=white] (6.8,2.5) circle (2 pt);
	\draw[ball color=black] (7.4,0.8) circle (2 pt); 
	\draw[black, thick,fill=Salmon!60] (7.4,2.2) circle (2 pt);
	\draw[ball color=black] (7.7,0.7) circle (2 pt); 
	\draw[black, thick,fill=Salmon!60] (7.7,2.1) circle (2 pt);
	\draw[ball color=black] (7.9,0.65) circle (2 pt); 
	\draw[black, thick,fill=Salmon!60] (7.9,2.1) circle (2 pt);
	\draw[ball color=black] (8.1,0.6) circle (2 pt); 
	\draw[black, thick,fill=Salmon!60] (8.1,2.15) circle (2 pt);
	\node[black] at (6,0.9){\small{$x(\tau_k)$}};
	\node[black] at (4.2,3){\small{$x(\tau_k^-)$}};
	\node[red] at (7.2,0.2){\small{$\setT$}};
	\node[blue!50!black] at (3.5,0){$\setX$};
	\end{tikzpicture}	
	}
	\caption{Control objective for impulsively controlled systems. The solid points represent the states at times $\tau_k$, $k \in \mathbb{N}$, just after the jump, while the empty points represents the states at times $\tau_k^-$, just before the jumps. The solid black lines represent the free response of the system, at each interval $t \in [\tau_k,\tau_{k+1})$, $k \in \mathbb{N}$.}
	\label{fig:EschPlotTar}
	\end{subfigure}
	\hfill
	\begin{subfigure}[b]{0.45\textwidth}
	\centerline{
	\begin{tikzpicture}[scale=1]
	\draw[blue!50!black, line width=1.0pt] (3,3.5) --  (9.3,3.5) -- (9.3,-0.5) -- (3,-0.5) -- cycle;
	\draw[red, line width=1.0pt, fill=Salmon!60] (7,-0.25) --  (6.3,1.2) -- (6.7,2) -- (8.5,2.8) -- (9.2,1) -- (9,0.5) -- cycle;
	\draw[black, line width=1.0pt] (5,1.5) .. controls (4.5,1.8) and (5.8,2.5) .. (4.8,3);
	\draw[black, line width=0.5pt, dashed] (4.8,3) -- (5.8,1.3); 
	\draw[black, line width=1.0pt] (5.8,1.3) .. controls (5.5,1.4) and (6.3,2.5) .. (6,2.9);
	\draw[black, line width=0.5pt, dashed] (6,2.9) -- (6.8,1); 
	\draw[black, line width=1.0pt] (6.8,1) .. controls (6.6,1.1) and (7.3,2) .. (6.8,2.5);
	\draw[black, line width=0.5pt, dashed] (6.8,2.5) -- (7.4,0.8); 
	\draw[black, line width=1.0pt] (7.4,0.8) .. controls (7.8,1.2) and (7.6,1.4) .. (7.4,2.2);
	\draw[black, line width=0.5pt, dashed] (7.4,2.2) -- (7.7,0.7); 
	\draw[black, line width=1.0pt] (7.7,0.7) .. controls (8.2,1.2) .. (7.7,2.1);
	\draw[black, line width=0.5pt, dashed] (7.7,2.1) -- (8.1,0.6); 
	\draw[black, line width=1.0pt] (8.1,0.6) .. controls  (8.8,1.2) .. (8.15,1.9);
	\draw[black, line width=0.5pt, dashed] (8.15,1.9) -- (8.5,0.6); 
	\draw[black, line width=1.0pt] (8.5,0.6) .. controls (9.3,1.2) .. (8.55,1.7);
	\draw[black, line width=0.5pt, dashed] (8.55,1.7) -- (8.8,0.65); 
	\draw[black, line width=1.0pt] (8.8,0.65) .. controls (9.7,1.1) .. (8.85,1.6);
	\draw[ball color=black] (5,1.5) circle (2 pt); 
	\draw[black, thick,fill=white] (4.8,3) circle (2 pt);
	\draw[ball color=black] (5.8,1.3) circle (2 pt); 
	\draw[black, thick,fill=white] (6,2.9) circle (2 pt);
	\draw[ball color=black] (6.8,1) circle (2 pt); 
	\draw[black, thick,fill=white] (6.8,2.5) circle (2 pt);
	\draw[ball color=black] (7.4,0.8) circle (2 pt); 
	\draw[black, thick,fill=Salmon!60] (7.4,2.2) circle (2 pt);
	\draw[ball color=black] (7.7,0.7) circle (2 pt); 
	\draw[black, thick,fill=Salmon!60] (7.7,2.1) circle (2 pt);
	\draw[ball color=black] (8.1,0.6) circle (2 pt); 
	\draw[black, thick,fill=Salmon!60] (8.15,1.9) circle (2 pt);
	\draw[ball color=black] (8.5,0.6) circle (2 pt); 
	\draw[black, thick,fill=Salmon!60] (8.55,1.7) circle (2 pt);
	\draw[ball color=black] (8.8,0.65) circle (2 pt); 
	\draw[black, thick,fill=Salmon!60] (8.85,1.6) circle (2 pt);
	\node[black] at (6,0.9){\small{$x(\tau_k)$}};
	\node[black] at (4.2,3){\small{$x(\tau_k^-)$}};
	\node[red] at (7.2,0.2){\small{$\setT$}};
	\node[blue!50!black] at (3.5,0){$\setX$};
	\end{tikzpicture}	
	}
	\caption{Inadmissible trajectory for impulsively controlled systems. Although the states sampled at times $\tau_k$, $k \in \mathbb{N}$, just after the jump, and at times $\tau_k^-$, just before the jumps are feasible and converge to the target zone, the state trajectories (solid black lines) leave the target and the feasible sets.}
	\label{fig:inadmissible}
	\end{subfigure}
	\caption{Illustration of the sets and trajectories for the zone tracking problem on impulsively controlled systems.}
	\label{fig:contprob}
\end{figure}

\begin{figure}[!t]
	\centerline{
	\begin{tikzpicture}[scale=1]
	\draw[blue!50!black, line width=1.0pt,fill=RoyalBlue!60] (3,3.5) --  (10,3.5) -- (10,-0.5) -- (3,-0.5) -- cycle;
	\draw[Goldenrod!50!black,line width=1.5pt,fill=Goldenrod] (3,-0.5) .. controls (3.3,3) ..  (10,3.5) .. controls (10,-0.5) .. (9,-0.5) --  cycle;
	\draw[red, line width=1.5pt, fill=Salmon!60] (6.6,-0.4) --  (6,1.2) -- (6.7,2.2) -- (8.8,2.8) -- (9.6,0.6) -- cycle;
    \draw[SeaGreen!50!black,line width=1.5pt,fill=SeaGreen] (6.65,-0.35) -- (6.7,0.43) .. controls (7,2.3) .. (9,1.5)  --  (8.43,0.31) .. controls (7.5,0) .. (6.68,-0.35);
    \draw[LimeGreen!50!black,line width=1.5pt,fill=LimeGreen] (6.7,0.43) .. controls (7.5,1.5) .. (9,1.5)  --  (8.43,0.31) .. controls (7.5,0) .. (6.7,0.43);
	\draw[BurntOrange,line width=1.5pt,dashed] (3.1,0.65)--(10,0.2);
	\draw[BurntOrange,line width=2pt] (6.7,0.43)--(8.43,0.31);
	\node[black] at (4.4,2.4){\small{$\setXA$}};
	\node[black] at (8.5,2.4){\small{$\setT$}};
	\node[black] at (7.5,1.7){\small{$\setFT$}};
	\node[black] at (8.4,1.2){\small{$\setTAinv$}};
	\node[black] at (7.5,0.7){\small{$\setTs$}};
	\node[black] at (3.5,3){$\setX$};
	\node[black] at (4.8,0.8){$\setXs$};
	\end{tikzpicture}
	}
	\caption{Illustrative plot of the sets of interest in the problem description: $\setX$ feasible set  (blue), $\setXA$ the admissible set of $\setX$ (yellow), $\setT$ the valid target set (red), $\setXs$ and $\setTs$ the ICES in  $\setX$ and $\setT$, respectively (the dash and solid orange lines), $\setFT$ the admissible set of target set (cyan) and the impulsive invariant set $\setTinv$ (green).}
	\label{fig:EschPlotOrig}
\end{figure}
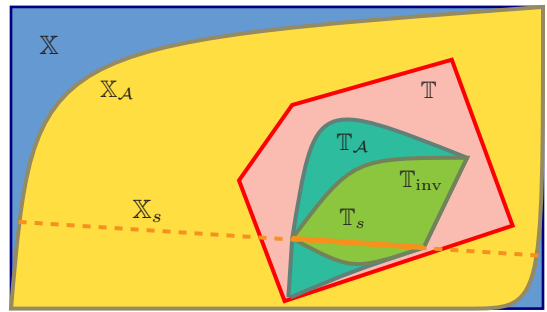

A direct consequence of Definition \ref{defi:targ} and Theorem \ref{theo:MainImp}, is stated next. 
\begin{cor}\label{cor:targequil}
	Consider the ICSys \eqref{ec:sysIm} and a valid target set $\setT$ as the one defined in Definition \ref{defi:targ}. 
	Then, $\setT$ contains a nonempty ICES, say $\setTs$.
\end{cor}
Figure \ref{fig:contprob} shows an schematic plot of the control problem for ICSys, while Figure \ref{fig:EschPlotOrig} describes all the equilibrium and invariant sets sets involved in it.

\section{Characterization of the controlled invariant sets for ICSys}

In this section it is shown how to characterize impulsive control invariant sets of ICSys \eqref{ec:sysIm}, by means of control invariant sets of its sampled version, DICSys \eqref{ec:sysDisc}. Next, a key theorem of the article is introduced:


\begin{theorem}\label{theo:Main1}
	Consider the ICSys \eqref{ec:sysIm} and a valid target set $\setT \subset \setX$. If a CIS $\setTinvD$ of the DICSys \eqref{ec:sysDisc} is a subset of $\setT_{\mathcal{A}}$, then it is also an ICIS of the ICSys \eqref{ec:sysIm}. We denote this set as $\setTinv$ for clarity.
\end{theorem}
\begin{proof}
	Given that $\setTinv \subset \setTA \subseteq \setT$, then $\{e^{At}x:~t \in [0,T]\} \in \setT$ for every $x \in \setTinv$. Particularly, $e^{AT}x$ is in $\setT$. Furthermore, since $\setTinv$ is a CIS for the DICSys \eqref{ec:sysDisc}, then it there exists some $u \in \setU$ such that $x^+:=e^{AT} x + Bu$, is in $\setTinv$, which concludes the proof.
\end{proof}
\begin{rem}
 Note that if a typical impulsive system representation - as the ones presented in \cite{Yang01,rivadeneira2018control,sopasakis2015model} - is used to describe the ICSys \eqref{ec:sysIm}, almost all the previous definitions are still valid, but Theorem \ref{theo:Main1} is no longer true. Indeed, if the free response is considered after the input jump, a CIS in $\setTA$ is not necessarily an ICIS. This is so because in that case it is not possible to characterize a single set $\{e^{At} \setTinv,~t \in [0,T]\}$ which is in $\setT$, i.e., CIS condition only says that for each $x \in \setTinv $ exists an $u \in \setU$ such that $e^{AT}(x+Bu) \in \setTinv$, but the inputs are in general different for different states. On the other side, with representation \eqref{ec:sysIm}, the input effect occurs after the free response, and the existence of an $u \in \setU$ that take the state back to $\setTinv$ is enough. 
\end{rem}
\begin{rem}
Note also that if a given set $\setYinv$ is a CIS but is not contained in $\setTA$, then the set $\setYinv \cap \setTA$ is not necessarily a CIS and, so, the hypothesis of Theorem \ref{theo:Main1} are no longer fulfilled (and $\setYinv$ is not an ICIS). In other words, for a set $\setYinv$ to be an ICIS, it should be computed from the beginning as a subset of $\setTA$, which is a point to be discussed in the next section.
\end{rem}
%
%

Figure \ref{fig:EschFig} shows an schematic plot of a CIS $\setTinv  \subset \setTA$, which is also an ICIS. 

\begin{figure}[!t]
	\centerline{\begin{tikzpicture}[scale=1.2]
	\draw[red, line width=1.5pt, fill=Salmon!60] (6,-0.5) --  (6,3) -- (7,4) -- (12,4) -- (12, 0.5)-- (11,-0.5) -- cycle;
    \draw[SeaGreen!50!black,line width=1.5pt,fill=SeaGreen] (6.65,-0.35) -- (6.7,0.43) .. controls (7,2.3) .. (9,1.5)  --  (8.43,0.31) .. controls (7.5,0) .. (6.68,-0.35);
    \draw[LimeGreen!50!black,line width=1.5pt,fill=LimeGreen] (6.7,0.43) .. controls (7.5,1.5) .. (9,1.5)  --  (8.43,0.31) .. controls (7.5,0) .. (6.7,0.43);
	\draw[BurntOrange,line width=1.5pt,dashed] (5.5,0.55)--(12.5,-0.05);
	\draw[BurntOrange,line width=2pt] (6.7,0.44)--(8.43,0.31);
    \draw[LimeGreen!50!black,line width=1.5pt,fill=LimeGreen] (9.2,2.43) .. controls (10,3.5) .. (11.5,3.5)  --  (10.93,2.31) .. controls (10,2) .. (9.2,2.43);
	\draw[BurntOrange,line width=2pt] (9.2,2.44)--(10.93,2.31);
 	\draw[->,black, line width=1.0pt] (7.8,0.37) .. controls (10,1.2) .. (10.3,2.31); 
 	\draw[->,black, line width=0.5pt, dashed] (10.25,2.32) -- (7.85,0.42); 
 	\draw[->,black, line width=1.0pt] (7.2,0.8) .. controls (8.5,3) .. (9.95,3); 
 	\draw[->,black, line width=0.5pt, dashed] (10,3) -- (7.85,1.25); 
 	\draw[->,black, line width=1.0pt] (9.5,0.23) .. controls (11,0.5) and (12,1) .. (12.28,2.23); 
 	\draw[->,black, line width=0.5pt, dashed] (12.25,2.25) -- (9.55,0.28); 
 	\draw[ball color=black] (7.8,0.37) circle (2 pt); 
 	\draw[black, thick] (10.3,2.37) circle (2 pt); 
 	\draw[ball color=black] (7.2,0.8) circle (2 pt); 
	\draw[black, thick] (10,3) circle (2 pt); 
	\draw[ball color=black] (7.8,1.2) circle (2 pt); 
	\draw[ball color=black] (9.5,0.23) circle (2 pt); 
 	\draw[black, thick] (12.3,2.3) circle (2 pt); 
	\node[black] at (6.8,3.2){$\setT$};
	\node[black] at (7.35,1.8){\small{$\setTA$}};
	\node[black] at (8.4,1.29){\small{$\setTinvD$}};
	\node[black] at (7.8,0.6){\small{$x_s$}};
	\node[black] at (10.3,2.6){\small{$e^{AT}x_s$}};
	\node[black] at (10.85,3.31){\small{$e^{AT}\setTinv$}};
	\node[black] at (7.25,0.58){\small{$\setTs$}};
	\node[black] at (5.5,0.8){$\setXs$};
	\end{tikzpicture}}
	\caption{Illustrative plot of . A CIS $\setTinvD$ contained in $\setTA$ is also an ICIS, $\setTinv$ (green).} 
	\label{fig:EschFig}
\end{figure}

 \subsection{Admissible Set of $\setT$ in spectrahedron representation} \label{sec:AdmCharact}

The idea now is to characterize the admissible set of $\setT$, $\setTA$, as the intersection of an uncountable number of polytopic sets (in the form of hyperplane constraints), which is known as spectrahedron representation.
Let $\setT$ be a polytopic set described by 
\begin{equation}
    \setT = \{x \in \setX : H x \leq v\},
    \label{setTpolytopic}
\end{equation}
with $H \in \mathbb{R}^{\ell \times n}$ and $v \in \mathbb{R}^{\ell \times 1}$, where $\ell$ is a minimal number of hyperplanes that describe $\setT$.
Next, the set $\setFT$ will be described in terms of linear matrix inequalities (LMI) as in \cite{henrion2005control}.
For the sake of completeness, the analysis is included in this work.
\begin{assumpt}
(i) The eigenvalues of matrix $A$ are rational numbers $\lambda_r \in \mathbb{Q}$, with no imaginary part, so that $\lambda_r = \frac{\eta_r}{\rho}$ where $\eta_r \in \setZ, \rho \in \N$ is the least common denominator of the eigenvalues and $r = 1, \cdots, n.$ The eigenvalues are ranked in increasing order, such that $\eta_1$ is the smallest (possibly negative) integer and $\eta_n$ is the  largest (possibly positive) integer.
(ii) The eigenvalues of A are all distinct.
\end{assumpt}

For $\setT$ described as in \eqref{setTpolytopic}, the admissible set is given by

\begin{equation}
    \setFT = \{x \in \setX: H\Phi(t)x \leq v, \forall t \in [0,T]\},
\label{eq:setFTpoly}
\end{equation}
where $\Phi(t) = e^{At}$ is the transition matrix.
Addressing inequality \eqref{eq:setFTpoly} row by row, it comes
\begin{equation}
    h_i \sum_{j=1}^{n}\Phi_j(t) x_j\leq v_i, ~ i = 1,\cdots,l
    \label{eq:setFTrows}
\end{equation}
where $h_i$ is the $i$-th row of $H$, $v_i$ is the $i$-th component of $v$ and $\Phi_j$ the $j$-th column of the transition matrix. Summation of the columns of the transition matrices scaled by the corresponding state component is just the column-wise expression of the matrix-by-vector multiplication. This can be interpreted as the contribution of each state component to the dynamic evolution through the transition $\Phi$. 
The transition matrix can be expressed in its modal form (see Appendix \ref{spectral-representation}) as 
\begin{equation*}
    \Phi(t) = \sum_{r=1}^{n}\phi^r e^{\lambda_r t},
\end{equation*}
where $\lambda_r$ is the $r$-th eigenvalue of the dynamic matrix $A$, and $\phi^r$ is the $r$-th matrix obtained from the modal decomposition of $\Phi$.
This can be used to rewrite \eqref{eq:setFTrows} as
\begin{equation}
    \sum_{r=1}^n
    \beta_r^i x e^{\lambda_r t} \leq v_i, ~ t \in [0,T], ~ i = 1,\cdots,l,
    \label{eq:setFTineq}
\end{equation}
with $\beta_{r}^i = h_{i} \phi^r$, since $h_i$ can be moved inside the sum in \eqref{eq:setFTrows}.

Proposing the change of variables $w = e^{-\frac{1}{\rho}t}$, follows $e^{\lambda_r t} = (e^{-\frac{1}{\rho}t})^{-\eta_r} = w^{-\eta_r}$, the interval on the new variable is $w \in [W,1]$ with $ W = e^{-\frac{T}{\rho}}$ and the inequalities \eqref{eq:setFTineq} can be expressed as
\begin{equation*}
    \sum_{r=1}^n \gamma^{i}_r(x)w^{-\eta_r}-v_i \leq 0, ~ w \in [W,1], ~ i = 1,\cdots,l,
    \end{equation*}
where $\gamma^{i}_r(x) = \beta_r^i x$
are linear functions of the initial state $x$ and $\eta_r \in \mathbb{Z}$, for $r \in 1, \cdots, n$.


At this point, we are interested in expressing the constraint as a polynomial with positive powers of the variable $w$. 
If $\eta_n$ ---which is the largest exponent--- is positive, we can multiply the previous expression by $w^{-\eta}$ and its reciprocal such that all the resulting terms in the sum correspond to non-negative integer powers of $w$.
We denote this \textit{degree-shift} as $\bar{\eta}:=\max\{0,\eta_n\}$, and it comes that

    \begin{equation*}
\frac{1}{w^{\bar{\eta}}}\left[\sum_{r=1}^n \gamma^{j}_r(x)w^{-\eta_r+\bar{\eta}}-v_i w^{\bar{\eta}}\right] \leq 0, ~ w \in [W,1], ~ i = 1,\cdots,l.
    \end{equation*}    
Noting that the factor $\frac{1}{w^{\bar{\eta}}}$ is non-negative in the interval $w \in [W,1]$, the previous condition can be written as a polynomial positivity condition as follows

\begin{equation}
P_i(w) = \sum_{d=0}^{\bar{\eta}}\pi_{i,d}(x) w^d \geq 0,  ~ w \in [W,1], ~ i = 1,\cdots,l,
\label{eq:inequalityP}
\end{equation}
where $\pi_{i,d}(x) = -\gamma_d^j(x) + \delta_{d,\bar{\eta}} v_i$, for $d \in \{0, \cdots, \eta_n+\bar{\eta}\}$, $\delta_{i,j}$ denotes the Kronecker's delta function (the case of null eigenvalues is considered in a single formulation). 
The $v_i$ coefficient is then considered for $d = \bar{\eta}$.
    
The Lukasz-Markov theorem states that a polynomial $P_i(w)$ is non-negative if and only if it can be written as a weighted sum of squares (see Appendix \ref{sec:SoS}).
Also, from \cite[Lemma 2]{Henrion_2005} it can be stated that inequality \eqref{eq:inequalityP} will be satisfied (i.e. the polynomial $P_i$ is non-negative) if its vector of coefficients $\pi_{i,d}(x)$
is the image of two positive semi-definite matrices $Y_{i,1}$ and $Y_{i,2}$ through linear operators $\Lambda^*_{i,1}$ and $\Lambda^*_{i,2}$,
as follows
    \begin{equation}
        P_i(x) = \Lambda^*_{i,1}(Y_{i,1}) + \Lambda^*_{i,2}(Y_{i,2}), ~ Y_{i,1},Y_{i,2} \succeq 0.
        \label{eq:lmicondition}
    \end{equation}
The operators $\Lambda^*_{i,1}$ and $\Lambda^*_{i,2}$ are defined in \cite{Nesterov2000}. 
Then, the set $\setFT$ is described by a semi-algebraic set where each element $x$ fulfills the condition \eqref{eq:lmicondition}:
    \begin{equation}
        \setFT = \{x \in \setX : P(x) = \Lambda^*_1(Y_1) + \Lambda^*_2(Y_2), ~ Y_1,Y_2 \succeq 0\}.
        \label{eq:setFT_LMI}
    \end{equation}
    
\subsection{ICES computation}
The ICES can be computed as the intersection of $\setFT$ and the CES of the DICSys system \eqref{ec:sysDisc}. The resulting ICES can be interpreted as a slice of the spectrahedric set $\setFT$ defined by the LMI \eqref{eq:setFT_LMI} and the CES defined by equality constraints.

\subsection{CIS computation - Polytopic Inner Approximation}
Tools for invariant set computation for generic convex sets (e.g. spectrahedra) are not widely available. 
Nevertheless, there are efficient algorithms for computing invariant sets with respect to polytopes. 
Note that a polytope with all its extreme points contained in a spectrahedron can be considered as an inner approximation. 
Then, inner polytopic approximation of the semidefinite constraints enable the use of classical algorithms for computation of invariant sets.\footnote{The subject of polytopic approximation of convex bodies is well studied and is beyond the scope of this work, but the reader is referred to \cite{gruber1983approximation,gruber1993aspects,bronstein2008approximation} for a thorough review on the subject.}

Let us denote $\setTAapprox$ a polytopic inner approximation of $\setFT$. By means of classical algorithms (\cite{MPT3,KerriganPHD00}), the computation the maximal CIS for the DICSys \eqref{ec:sysDisc} contained in $\setTAapprox$ is an invariant set for the DICSys. 
Also, note that this set is also an inner approximation to the maximal ICIS for the ICSys \eqref{ec:sysIm} w.r.t. $\setT$, and the continuous time trajectories lie within $\setT$ at all times.

The set is typically computed recursively, where the successive controllable sets are obtained by the following recursion
	\begin{equation*} 
	S_{k+1} = \{ x  |  A^d x + B^d u \in S_k, u \in \setU \},
	\end{equation*}
and terminating once $S_{k+1}=S_k$. 

In order to obtain the invariant set within the target zone $\setT$, the algorithm is initialized with $S_0 = \setTAapprox$.  For the computation of the feasible set  with respect to $\setX$, the algorithm is initialized  with $S_0 = \setFXapprox$, instead. 
The resulting ICISs (which are also CISs, as previously explained) contained in $\setX$ and $\setT$ are denoted by $\setPolyXinv$ and $\setPolyXinvT$ respectively.

The next section takes advantage of it to design a model predictive controller.

\section{Zone MPC control}
This section is devoted to introduce zone model predictive control (zMPC) formulations that make an explicit use of the ICIS characterized in the previous sections. The control objective --- as stated in the Control Problem definition of Subsection \ref{sec:contrprob} --- is to feasibly steer the impulsively controlled system \eqref{ec:sysIm} to the target zone $\setT$ (which does not necessarily contains the origin) and keep it there indefinitely.

Two different controllers are presented next: the first one is an MPC based on the use of artificial equilibrium variables \cite{FerramoscaJPC2012}, which ensures asymptotic stability. The second one is a set-based MPC \cite{AndersonOCAM18,liu2019zone}, which guarantees, in addition, finite-time convergence.

\subsection{Zone MPC based on Artificial Variables}
The DICSys \eqref{eq:dsys} is used for predictions and an invariant set $\setTinvD$ contained in $\setTA$ (prferably the largest) acts as target set.
The cost function to be minimized on-line by the MPC is given by
%
 \begin{eqnarray}\label{MPC_cost}
 J_N(x,\mathbf{u},x_s,u_s) &=& \sum_{j=0}^{N-1} \|x(\tau_j) - x_s\|^Q_2 + \|u(\tau_j)-u_s\|^R_2 + \gamma \mathrm{dist}_{\setFT}(x_s) 
 \end{eqnarray}
where $x = {x}(\tau_0)$ represents the current state, $\mathbf{u}=\{u(\tau_0),u(\tau_1),\dots,u(\tau_{N-1})\}$
is the predicted sequence of inputs, 
and $Q$ and $R$ are positive definite and semi-definite matrices, respectively. 

\begin{rem}\label{rem:distfunc}
The term $ \mathrm{dist}_{\setFT}(x_s)$ is approximated by using an additional optimization variable $x^*$, which is forced to be in $\setFT$, and by including the cost term $\|x-x^*\|_2$. It is computed as the optimization problem given by $\mathrm{dist}_{\setFT}(x_s)=\min_{x^*\in\setFT}\|x_s-x^*\|$.
	The problem in included into the MPC and is solved simultaneously by modifying the cost function, which now reads 
		 \begin{equation}\label{MPC_cost1}
		 \begin{aligned}
		\displaystyle J_N(x,\mathbf{u},x_s,u_s,x^*) & = \sum_{j=0}^{N-1} \|x(\tau_j) - x_s\|^Q_2 + \|u(\tau_j)-u_s\|^R_2 \\ &+ \|x_s-x^*\|^{Q_O}_2,
		 \end{aligned}
		 \end{equation}
		and constraint $x^* \in \setFT$ is included in the optimization \eqref{MPC_problem},
	introduced next.
\end{rem}

The optimization problem to be solved at each time $\tau_k$ is as follows:
\begin{subequations}\label{MPC_problem}
\begin{align}
	\min_{\mathbf{u},x_s,u_s,x^*}  & J_N(x,\mathbf{u},x_s,u_s,x^*)\\
	s.t. & \nonumber \\
	& x(\tau_0)=x\\
	& x(\tau_{j+1})= A^d x(\tau_j)+B^d u(\tau_j), &j \in 0, \cdots, {N-1}\\
	& x(\tau_{j+1}) \in \setXA, & j \in 0, \cdots, {N-1}\\
	& u(\tau_j) \in \setU, & j \in 0, \cdots, {N-1}\\
	& A^d x_s + B^d u_s = x_s, ~ x_s \in \setXs, \label{eq:invariance.cons}\\
	& x(\tau_N)  = x_s, \label{eq:terminal.const}\\
	&x^* \in \setTA, \label{eq:dist.const}
	\end{align}
\end{subequations}

In the latter optimization problem, $x$ is the optimization parameter while $\mathbf{u},x_s,u_s,x^*$ are the optimization variables.
Constraint \ref{eq:invariance.cons} forces the pair of additional variables $(x_s,u_s)$ to be contained in (an inner polytopic approximation of) the lifted admissible invariant set of the sampled system.
The terminal constraint \ref{eq:terminal.const} forces the last state on the control horizon to reach the equilibrium set, as this will be used to ensure stability.
Constraint \ref{eq:dist.const} is included to approximate the distance function in the cost by a quadratic term (Remark \ref{rem:distfunc}).
%
%
%

Once the MPC problem is solved at time $\tau_k$, the optimal solution is given by the optimal input sequence 
\begin{equation}
\mathbf{u^0}(x) = \{u^0(x,\tau_0),u^0(x,\tau_1),\dots,u^0(x,\tau_{N-1})\},    
\end{equation}
and the optimal additional variables $(x_s^0(x),u_s^0(x))$,
while the optimal cost is denoted as $J_N^0(x) := J_N(x,\mathbf{u}^0,x_s^0,u_s^0)$.
The control law, derived from the application of a receding horizon control policy (RHC),
is given by $u(\tau_k) = \kappa_{\mathrm{MPC}}(x(\tau_k)) = u^0(x,\tau_0)$, where $u^0(x,\tau_0)$ is the first control action in $\mathbf{u^0}(x)$. 

The stabilizing properties of the controller are summarized in the following Property.
	\begin{propt}
     The IES $\setTs$ is asymptotically stable for the closed-loop system ICS, 
     \begin{equation}\label{ec:closloop}
     \begin{cases}
     \dot{x}(t)  & =  A x(t),~~~ t \neq \tau_k, \\ 
     x(\tau_k^+) & =  x(\tau_k) +B \kappa_{\mathrm{MPC}}(x(\tau_k)),~~~~~ k \in \mathbb{N},
     \end{cases}
     \end{equation}
     with $x(0)=x_0$
     , and a domain of attraction given by $\setC_N(\setXs)$
     .
	\end{propt}
\begin{spf}
	The proof follows the steps of the so called MPC for tracking zone regions, \cite{FerraLi10,rivadeneira2018control}.
	The time of the closed-loop is denoted by $\tau_k$, as in \eqref{ec:closloop}, while
	the time for predictions, inside each optimization problem, is denoted by $\tau_j$.
	\\
	The recursive feasibility of the sequence of optimization problems
	follows from the fact that, for every $x \in \setC_N(\setXs)$, the terminal constraint forces the system $x(\tau_{j+1})=\phi(T)x(\tau_j)+Bu(\tau_j)$
	to reach an invariant set (for instance the equilibrium set, $\setXs^0$),
	at the end of the control horizon. So, if the solution of the optimization problem for
	$x$, at time $\tau_k$, is given by $\mathbf{u^0}(x)$, $x_s^0(x)$ and $u_s^0(x)$, then a feasible solution for 
	the state $x^+$, at time $\tau_{k+1}$, can be computed as 
	$\tilde{\mathbf{u}}(x^+)=\{u^0(x;\tau_1),u^0(\tau_2),\dots,u^0(\tau_{N-2}),\tilde{u}_s(x^+)\}$,
	$\tilde{x}_s(x^+) = x_s(x)$ and $\tilde{u}_s(x^+)=u_s(x)$.
	This feasible solution produces a feasible sequence of states, given by
	$\tilde{\mathbf{x}}(x^+)=\{x^0(x,\tau_1),$ $x^0(x,\tau_2),\dots,\tilde{x}_s(x^+),\tilde{x}_s(x^+)\}$.
	\\
	The attractivity of $\setXs$ follows from the fact that, $\tilde J_N(x^+) \leq J_N^0(x) - \alpha \|x - x_s\|_2 -\beta \|u-u_s\|_2$,
	where $\tilde J_N(x^+) = J_N(\tilde{\mathbf{u}}(x^+),\tilde{x}_s(x^+),\tilde{u}_s(x^+))$ and $u$ is the input injected to the system at time $\tau_k$. 
	Then, by optimality, it is $J_N^0(x^+) \leq \tilde J_N(x^+)$, which means that
	$J_N^0(\cdot)$ is a strictly decreasing positive function - i.e., $J_N^0(x^+) \leq  J_N^0(x)$ - that only stops to decrease	if $x=x(\tau_k)=x_s$ and $u(\tau_k)=u_s$.
	Furthermore, the fact that $x(\tau_k) \rightarrow x_s$	and $u(\tau_k) \rightarrow u_s$, as $k \rightarrow \infty$, implies that $x(\tau_k)$ tends also to $\setT_{\setZ}$ (by the effect of the cost term $\mathrm{dist}_{\setT_{\setZ}}(x_s)$, as stated in Lemmas 1, 2 and 3, in \cite{rivadeneira2018control}). 
	This way, $x(\tau_k)$ tends to the intersections of  $\setXs^{\circ}$ and $\setT_{\setZ}$, which represents the IES $\setXs$.
\end{spf}

\begin{rem}
	Note that it is not necessary to express explicitly the intersection to formulate the MPC. In fact, such intersection is implicit in the controller formulation, by means of the additional variables $(x_s,u_s)$ - that are forced
	to be in $\setXs^{\circ}$ -, and the cost term $\mathrm{dist}_{\setT_\setZ}(x_s)$ which steers the states after the discontinuities,
	$x^{\circ}(\tau_k)$, to $\setT_{\setZ}$.
\end{rem}

\begin{rem}
	In \cite{rivadeneira2018control} a target set $\setT$, already accounting for the properties of $\setXs$ w.r.t. $\setZ$
		(i.e., accounting for $\{\Phi(t)x_s,~t \in [0,T]\} \in \setZ$), needs to be outer-approximated by a polyhedron, and then explicitly 
		used in the controller formulation. In contrast, the proposed MPC steers the system to the exact set $\setXs$, without the need
	of explicitly compute it.
\end{rem}

\begin{rem}
Another benefits of the proposed MPC to be emphasized is that it steers the system to an equilibrium region that fulfill continuous-time constraints by only considering a sampled discrete-time system, as it is $x^\circ(\tau_{j+1})=\phi(T)x^\circ(\tau_j)+Bu(\tau_j)$.
\end{rem}

\subsection{Set-based Zone MPC For Impulsive Systems}
This formulation is based on the ones reported in \cite{AndersonOCAM18,liu2019zone}, which provide better transient performance of the zone tracking scheme, as they tracks an invariant set within the target zone rather than steady-state setpoints and, furthermore, thy may guarantee finite time convergence.
The objective is to minimize the distance between the predicted state and inputs trajectory with respect to the lifted set $\setZ$ of the states and their corresponding inputs, such that $\setZ = {(x,u) | x \in \setTinv, u \in \Ucal(x)}$.
The cost function to be minimized on-line by the set-based MPC is given by:
\begin{equation*}
    L_N = \sum_{j*1}^n \dist_\setZ((x(\tau_j),u(\tau_j)),
\end{equation*}
where the distance function is implemented as
\begin{equation*}
    \dist_\setZ(x,u) = \min_{(\mathbf{x}^*,\mathbf{u}^*) \in \setZ}\|x-x^*\|+\|u-u^*\|.
\end{equation*}
The proposed formulation is
\begin{subequations}\label{MPC_problem_2}
\begin{align}
	\min_{\mathbf{u},\mathbf{x^*},\mathbf{u^*}}  & L_N(x,\mathbf{u},\mathbf{x}^*,\mathbf{u}^*)\\
	s.t. & \nonumber \\
	& x(\tau_0)=x\\
	& x(\tau_{j+1})= A^d x(\tau_j)+B^d u(\tau_j), &j \in 0, \cdots, {N-1}\\
	& x(\tau_{j+1}) \in \setXA, & j \in 0, \cdots, {N-1}\\
	& u(\tau_j) \in \setU, & j \in 0, \cdots, {N-1}\\
	&x^*(\tau_j) \in \setTinv, & j \in 0, \cdots, {N-1} \label{eq:dist.const2}\\
	& u^*(\tau_j) \in \setU, & j \in 0, \cdots, {N-1}\\
	& A^d x^*(\tau_j) + B^d u^*(\tau_j) \in \setTinv, &j \in 0, \cdots, {N-1} \label{eq:invariance.cons2}\\
	& x(\tau_N)  = x^*(\tau_N). \label{eq:terminal.const2}
	\end{align}
\end{subequations}
Denoting the resulting receding horizon feedback law by $\kappa^*_{\mathrm{MPC}}$, the following property is obtained.
	\begin{propt}
     The IES $\setTs$ , is asymptotically stable for the closed-loop system ICS, 
     \begin{equation}\label{ec:closloop}
     \begin{cases}
     \dot{x}(t)  & =  A x(t),~~~ t \neq \tau_k, \\ 
     x(\tau_k^+) & =  x(\tau_k) +B \kappa^*_{\mathrm{MPC}}(x(\tau_k)),~~~~~ k \in \mathbb{N},
     \end{cases}
     \end{equation}
     with $x(0)=x_0$
     , and a domain of attraction given by $\setC_N(\setTinv)$
     .
	\end{propt}

\section{Example}
\begin{figure}
     \centering
     \begin{subfigure}[b]{0.3\textwidth}
         \centering
         \includegraphics[width=\textwidth]{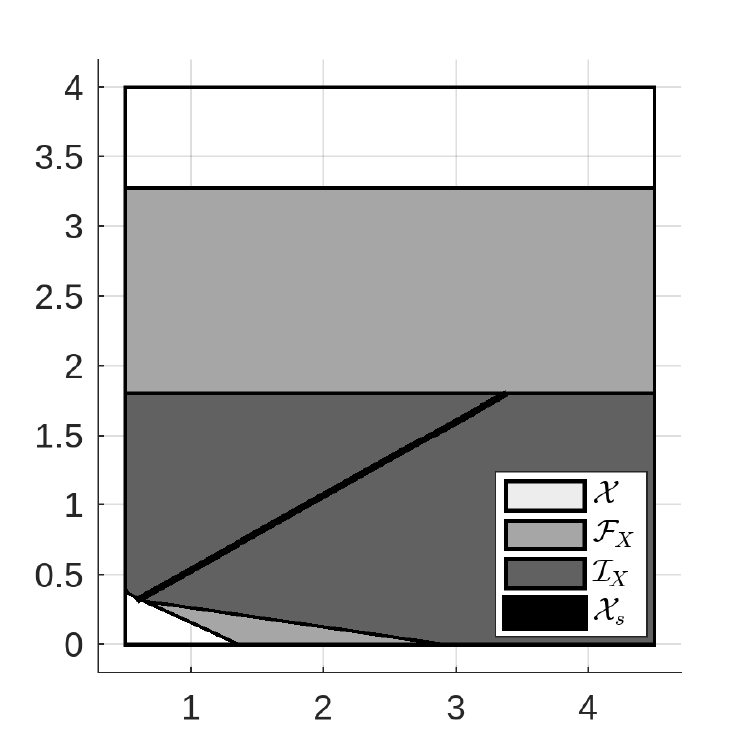}
         \caption{Plot of the feasible (white), Admissible (light gray) and Invariant with respect to the admissible (dark gray), admissible equilibrium (black), sets.}
         \label{fig:feasible-sets}
     \end{subfigure}
     \hfill
     \begin{subfigure}[b]{0.3\textwidth}
         \centering
         \includegraphics[width=\textwidth]{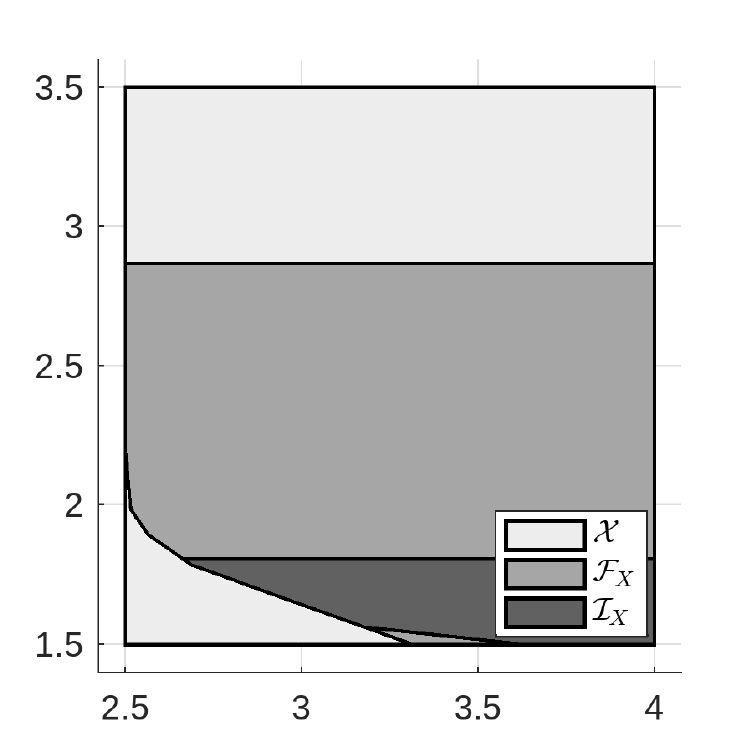}
         \caption{Plot of the target (light gray), Admissible (medium gray) and Invariant with respect to the admissible (dark gray), sets.}
         \label{fig:target-sets}
     \end{subfigure}
     \hfill
     \begin{subfigure}[b]{0.3\textwidth}
         \centering
         \includegraphics[width=\textwidth]{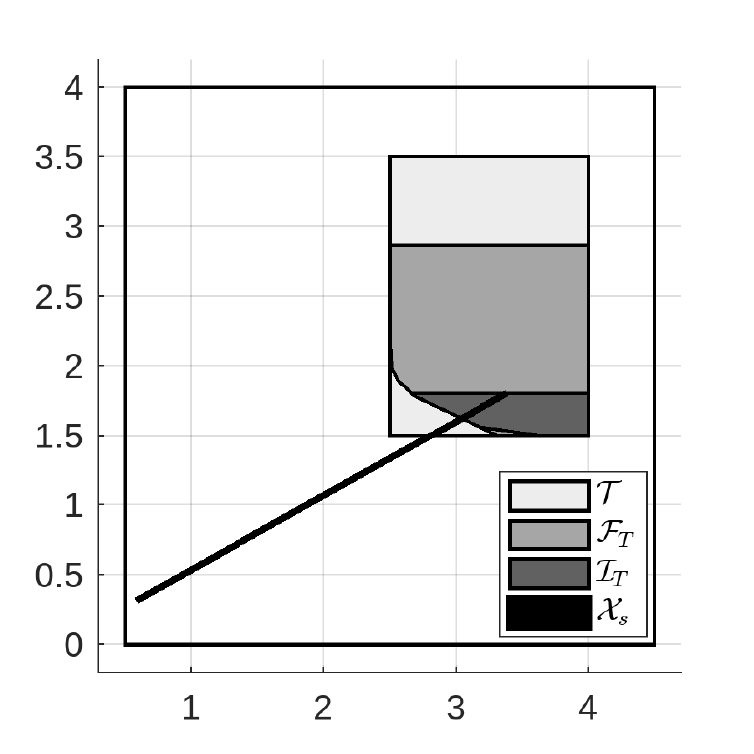}
         \caption{Target window and admissible equilibrium set.}
         \label{fig:target-and-equilibrium}
     \end{subfigure}
     \end{figure}

\begin{figure}
     \centering
     \begin{subfigure}[b]{0.45\textwidth}
        \centering
        \includegraphics[width=\textwidth]{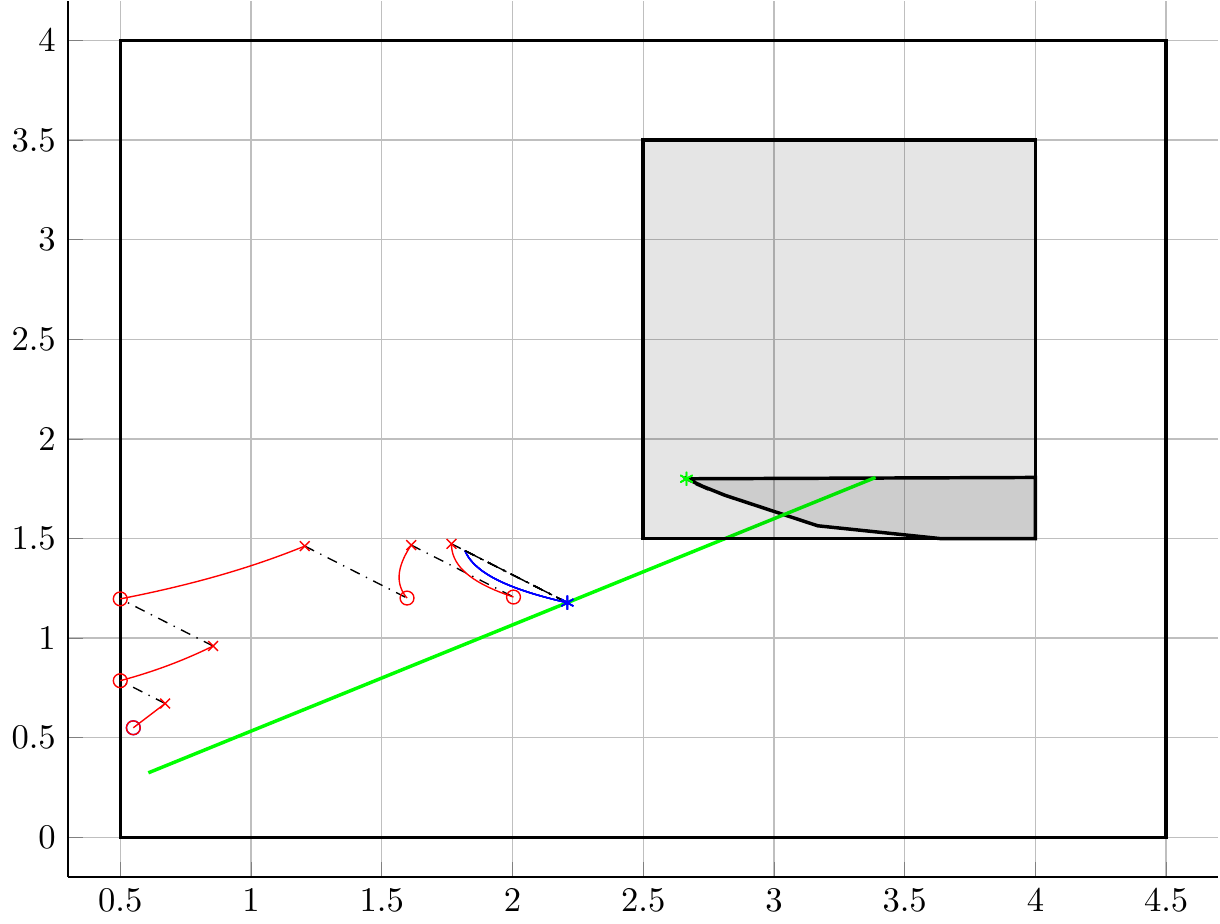}
        \caption{Optimal trajectory for initial state $(0.55,0.55)$. The red curves indicate the predicted continuous system evolution, the dotted lines the jumps on the state due to the impulsive inputs. The blue trajectory indicates the impulsive periodic trajectory and the green asterisk indicates the optimal $x^*$ for this solution.}
        \label{fig:plot-prediction}
    \end{subfigure}
         \hfill
     \begin{subfigure}[b]{0.45\textwidth}
        \centering
        \includegraphics[width=\textwidth]{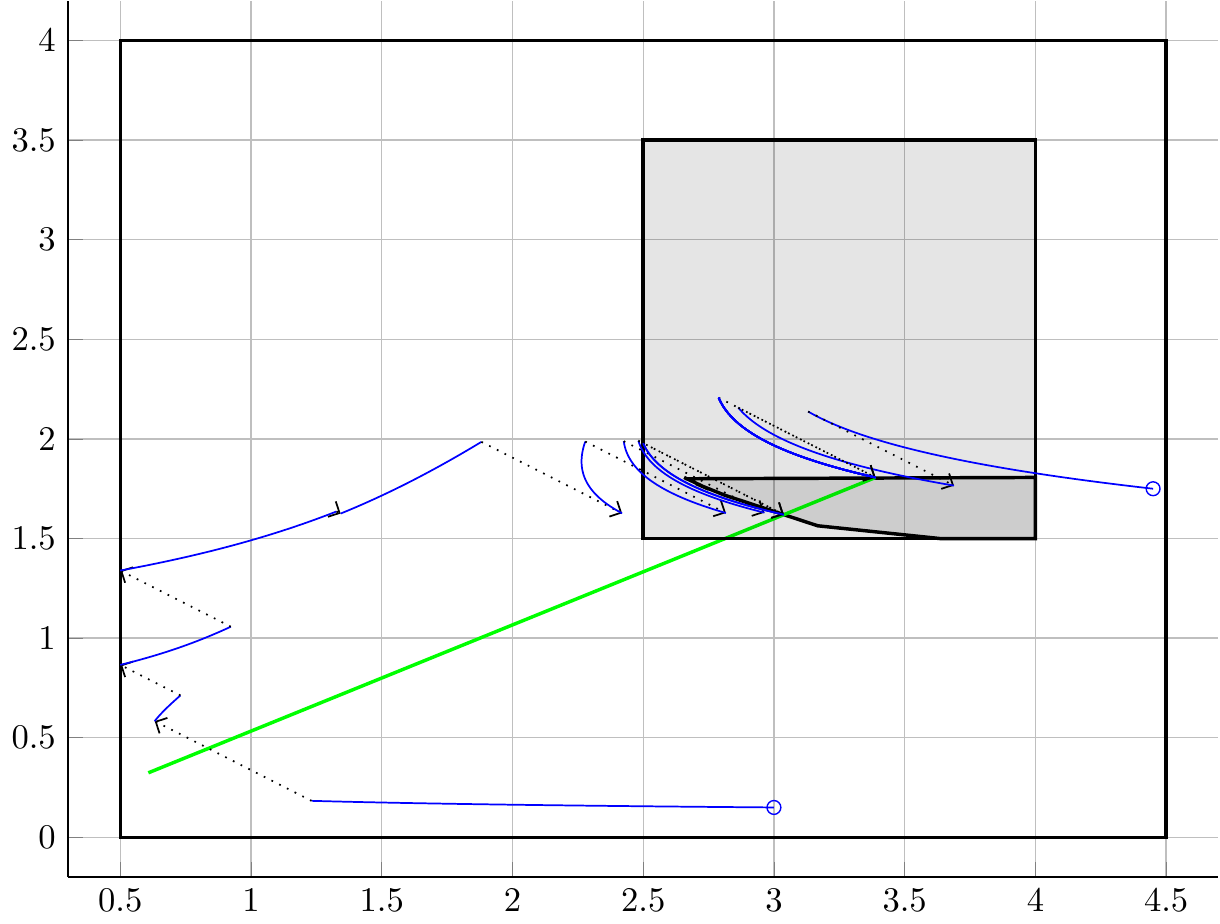}
        \caption{Simulation of 10 iterations of the closed-loop system, starting at two different initial states (3.0,0.15) and (4.45,1.75). The system converges in both cases to the zone, at different equilibrium states that correspond to each initial point.}
        \label{fig:plot-solutions}
    \end{subfigure}
    \caption{Single iteration and solution trajectories for the closed loop system. The target window is shown in light gray and the admissible invariant set within the set is shown in darker gray. Note that the continuous time trajectories fulfill the constraints at all times.}
    \end{figure}

Consider the linear system \ref{ec:sysIm} with dynamic and input matrices given by 
$$
    A = \begin{bmatrix}
    -1 & 1.2 \\ 0 & 0.2
    \end{bmatrix}, ~
    B = \begin{bmatrix}
    3 \\ -2
    \end{bmatrix}.$$
    The state is constrained to a box defined by  $0.5 \leq x_1 \leq 4.5$ and $0 \leq x_2 \leq 4$ and the input is bounded $u \in [-0.2, 0.2]$.
    The input impulses occur at intervals of $T = 1 [s]$.
    The target window is a box with $2.5 \leq x_1 \leq 4$ and $1.5 \leq x_2 \leq 3.5$.
    The admissible polytopic approximation to the admissible spectrahedron and invariant set with respect to the feasible state and target window are shown in figures \ref{fig:feasible-sets} and \ref{fig:target-sets}, respectively. The equilibrium set and target window are shown in \ref{fig:target-and-equilibrium}, and makes evident that the target window is a valid target set. 
    The eigenvalues are $-1$ and $.2$. The autonomous system is unstable, diverging naturally from the target set.
    The controller is defined with a horizon $N = 5 [s]$, and the weighting matrices are chosen as $Q = 1, R = 10, Q_f = 1, Q_O = 10$. The simulations are conducted using YALMIP \cite{lofberg2004yalmip}, MPT3 \cite{MPT3} and Mosek \cite{andersen2000mosek}. \\
    For illustrative purposes, the resulting optimal trajectories and values for auxiliary ingredients of the proposed controller with the system starting at $(0.55,0.55)$ are shown in Figure \ref{fig:plot-prediction}. The state trajectory is shown to evolve on a hybrid trajectory, continuously in the interval $(t_k, t_{k+1})$ and discrete jumps occur at $t_k$. 
    Also the state reaches an equilibrium point  of the discrete-time associated system given by $(x_s,u_s)$ shown in blue, which is a periodic solution in continuous time. 
    The distance function is implemented as indicated in (\ref{rem:distfunc}) and the resulting value of $x^*$ is shown in green, which is a state contained in the admissible invariant set contained in the target window. \\
    The resulting trajectories for the system starting at two different initial states $(3.0,0.15)$ and $(4.45,1.75)$ are shown in Figure \ref{fig:plot-solutions}. 
    They correspond to 10 iterations of the optimal controller implementation. Both trajectories converge to the target zone, more particularly to the admissible invariant set within the target. 
    In fact, both reach their respective admissible equilibrium solution, as expected. 
    An important remark is that the the steady-state solution for the discrete time systems -and the corresponding continuous time orbits- are not necessarily coincident, although they fulfill the control problem objectives.
    It is also worth noting that the initial points are contained in, and the sampled trajectory evolves within the invariant admissible set.

\section{Conclusions}
In this work, a exact characterization of admissible states for linear impulsively controlled systems is discussed, with application to systems with different rational real eigenvalues.
Using a description of the impulsive system where the impulsive input is applied at the end of each sampling period, the admissible states can be determined using a strategy based on the Lukasz-Markov theorem, exploiting and semidefinite  descriptions and sum-of-squares expressions of the system dynamics and constraints. 
This problem can be tractably solved using SDP and sparse computational techniques.
Conditions on the validity of a target zone are discussed, as well as a proof based on the Brouwer fixed-points theorem is provided. This demonstrates the existence of equilibrium points within an invariant set for linear systems, conforming a new approach to the proof using simple techniques. The strategies already available in the literature relied on the more general Kakutani fixed-point theorem for correspondences.

It is shown that the inclusion of an equilibrium set in the target zone implies the validity of the target set and follows that the system is able to be steered in such a way that its state remains inside the target zone.
Classic algorithms for the computation of an invariant set within a polytopic set is used, through inner approximation of the constrained admissible set, which is defined as spectrahedron, a type of convex set described by semidefinite inequalities. 
The resulting invariant set is used in a MPC formulation for zone tracking, and simulations are included to illustrate the effectiveness of the proposal. 
\section{Appendix A. Fixed point theorem and Berge's principle} \label{sec:AppA}

First, the well known Brouwer fixed-point is established.
\begin{theorem}[Fixed-point theorem, \cite{border1985fixed}] \label{theo:fixPoint}
	Given a compact and convex set $\setY \subset \mathbb{R}^n$ and a continuous function $f: \setY \rightarrow \setY$, then
	there is a point $x_{0} \in \setY$ such that $f(x_{0})=x_{0}$; that is, $x_0$ is a fixed point of $f$.
\end{theorem}
This theorem states that for any continuous function $f$ mapping a compact convex set to 
itself, there is at least one fixed point. For instance, every compact and convex invariant set for the autonomous 
system $x^+=A^d x$, say $\setXinv \subset \setX$, includes an equilibrium point, $x_s=A^d x_s$ 
(as stated in \cite{Blanchinibook15}, page 113, for both, the continuous and discrete-time cases.).

Next, the Berge's principle is established:
\begin{theorem}[Berge's principle \cite{berge1959}]\label{theo:Berge}
	Consider a variable $y \in \setY \subset \R^n$, a parameter $p \in \setP \subset \R^m$ (with $\setY$ and $\setP$
	compact and convex), and 
	the correspondence $\Gamma:\setP \rightrightarrows \setY$, defined by
	\begin{eqnarray}\label{ec:Gamma}
	\Gamma(p) := \argmin_{y \in \Phi(p)} J(y,p),~~~~ \forall ~p \in \setP,
	\end{eqnarray}
	where $\Phi(p) \subset \setY$, for all $p \in \setP$ ($\Phi$ is a correspondence assigning a set to each $p$, 
	i.e., $\Phi: \setP \rightrightarrows \setY$) and\\ 
	(i) $J(y,p)$ is a continuous function on $\setY$, for all $p \in \setP$, \\
	(ii) $\Phi(p)$ is a compact-valued correspondence (i.e., compact for each $p \in \setP$), and\\
	(iii) $\Phi(p)$ is continuous (i.e., both, upper and lower semicontinuous), on $\setP$.\\
	Then $\Gamma(p)$ is nonempty-valued, compact-valued and upper semicontinuous in $p \in \setP$.
\end{theorem}
Note that, if $\Gamma(p)$ is not only a correspondence, but also a function, the latter Theorem directly states that
$\Gamma(p)$ is a continuous function on $\setP$ (i.e., a function that is upper semicontinuous is continuous).
From the mathematical programming theory, conditions under which $\Gamma(p)$ is a single-valued correspondence (i.e., a function) are
the conditions under which, for each $p\in \setP$, there is a unique solution $\Gamma(p)$, i.e.:
\\
(i) $J(y,p)$ is a strictly quasi-convex function on $\setY$, for all $p \in \setP$, and
\\
(ii) $\Phi(p)$ is a nonempty-valued, convex-valued and compact-valued correspondence (i.e., nonempty, convex and compact, respectively, 
for each $p \in \setP$).
\\
So, by adding the latter conditions to the ones in Theorem \ref{theo:Berge}, we have the following
variant of the Berge's principle:
\begin{theorem}[Variant of Berge's principle]\label{theo:VarBerge}
	Consider the correspondence $\Gamma(p)$, as the one defined in \eqref{ec:Gamma}, with\\ 
	(i) $J(y,p)$ being a strictly quasi-convex, continuous function on $\setY$, for all $p \in \setP$, \\
	(ii) $\Phi(p)$ being a nonempty-valued, convex-valued and compact-valued correspondence on $\setP$, and\\
	(iii) $\Phi(p)$ being continuous (i.e., both, upper and lower semicontinuous) on $\setP$.\\
	Then $\Gamma(p)$ is a continuous function on $\setP$.
\end{theorem}
\section{Appendix B. Technical Lemmas} \label{sec:AppB}
First, consider the definitions of upper and lower semicontinuity (\cite{aubin2009set}): 
\begin{defi}[Upper semicontinuity]\label{defi:upsemi}
	A correspondence $c : \setX \rightrightarrows \setY$ is upper semicontinuous at $x_0$ if (and only if)
	for every $\epsilon >0$ it there exists $\delta>0$ such that:
	$x \in \setB_{\delta}(x_0)$ $\Rightarrow$ $c(x) \subseteq \setB_{\epsilon}(c(x_0))$.
\end{defi}
\begin{defi}[Lower semicontinuity]\label{defi:lowsemi}
	A correspondence $c : \setX \rightrightarrows \setY$ is lower semicontinuous at $x_0$ if (and only if)
	for every $\epsilon >0$ fulfilling $\setB_{\epsilon}(c(x)) \cap c(x_0) \neq \emptyset$ 
	it there exists $\delta>0$ such that for every $x \in \setB_{\delta}(x_0)$ it follows that
	$c(x) \cap \setB_{\epsilon}(c(x_0)) \neq \emptyset$.
\end{defi}

A useful property relating upper continuity and the graph of a correspondence is stated next.
\begin{lem}\label{lem:uppcont}
	Let $c : \setX \rightrightarrows \setY$ be a correspondence and let the image set ($c_{\setX}:= \cup_{x\in \setX} c(x)$) compact.
	Then, $c$ is upper semicontinuous on $\setX$ if and only if its graph ($Gr(c)$) is closed.	
\end{lem}

Next Lemmas provide useful properties of the correspondence $\setU(x) := \{ u \in \setU : A^dx + B^du \in \setXinv\}$ 
defined in \eqref{ec:u_of_x} (i.e., the set of all $u$'s that keep a particular $x \in \setXinv$
in $\setXinv$).

\begin{lem}\label{lem:imagclosure}
	Consider DICSys \eqref{ec:sysDisc} and a convex, compact, proper CIS, $\setXinv \subset X$.
	Let $\setU_{\setXinv} := \cup_{x \in \setXinv} \setU(x)$ be the image set of $\setU$.
	Then, $\setU_{\setXinv}$ is compact. 
\end{lem}

\begin{proof}
	Compactness means closeness and boundness. Let us recall that a set is closed if every convergent
	sequence in it converges to a point that is also in it.\\
	Let $\{u_k\}$ be a convergent input sequence in $\setU_{\setXinv}$, and let $u$ its limits (i.e., $\{u_k\} \rightarrow u$).
	Since $u_k \in \setU_{\setXinv}$, for all $k \in \mathbb{I}_{\infty}$, it there exist states $x_k$ 
	such that $u_k \in \setU(x_k)$. This means that $\{x_{k}\} \subset \setXinv$ and, given that
	$\setXinv$ is compact (by hypothesis), then it there exists a subsequence $\{x_{k_i}\}$ which
	converges to $x$, and $x \in \setXinv$.
	The fact that $u_{k_i} \in \setU(x_{k_i})$, for all $i \in \mathbb{I}_{\infty}$, means that $w_{k_i} := A x_{k_i} + B u_{k_i}$ is in
	$\setXinv$, for all $i \in \mathbb{I}_{\infty}$. Then, by continuity of the linear system, $\{A x_{k_i} + B u_{k_i}\} \rightarrow Ax+Bu$; that is,
	$\{w_{k_i}\}$ converges to $w:=Ax+Bu$.
	But $\{w_{k_i}\}$ is a convergent sequence in the closed set $\setXinv$, then its limits $w$ is also in $\setXinv$. Then, since $x \in \setXinv$ and $w \in \setXinv$, we have that $u \in \setU(x)$, which means that $u \in \setU_{\setXinv}$
	and so $\setU_{\setXinv}$ is closed. Given that $\setU_{\setXinv} \subset \setU$, and $\setU$ is bounded (compact), then $\setU_{\setXinv}$ is bounded. Finally, $\setU_{\setXinv}$ is both, closed and bounded, i.e., it is compact.
\end{proof}

\begin{lem}\label{lem:graphclosure}
	Consider DICSys \eqref{ec:sysDisc} and a convex, compact, proper CIS, $\setX \subset X$.
	Let $Gr(\setU) := \{(x,u)~:~ u \in \setU(x)\}$ be the graph of $\setU$.
	Then, $Gr(\setU)$ is closed and convex.
\end{lem}

\begin{proof}
	Let $\{(x_k,u_k)\}$ be a convergent sequence in $Gr(\setU)$, and let $(x,u)$ its limit (i.e., $\{(x_k,u_k)\} \rightarrow (x,u)$).
	This implies that both, $\{x_k\} \rightarrow x$ and $\{u_k\} \rightarrow u$.
	Furthermore, by the definition of the graph, $u_k \in \setU(x_k)$ for all $k \in \mathbb{I}_{\infty}$, which implies that
	$w_k:=Ax_k+Bu_k$ is in $\setX$. As $\{w_k\}$ converges to $w:=Ax+Bu$, and every convergent sequence in the closed set $\setX$
	converges to a point in $\setX$, then $w \in \setX$. Since $x\in\setX$ and $w\in\setX$, then $u \in \setU(x)$,
	and $(x,u) \in Gr(\setU)$.

	In order to prove that $Gr(\setU)$ is convex, take two different points $(x_1,u_1)$ and $(x_2,u_2)$ both in $Gr(\setU)$. Then both $Ax_1+Bu_1$ and $Ax_2+Bu_2$  belongs to $\setXinv$. Since $\setXinv$ is convex $\lambda(Ax_1+Bu_1) +(1-\lambda)(Ax_2+Bu_2)$ also belongs to $\setXinv$, for any $\lambda\in[0,1]$, which means that $A(\lambda x_1 +(1-\lambda)x_2)+B(\lambda u_1 +(1-\lambda)u_2) \in \setXinv$. Again since  $\setXinv$ is convex $\lambda x_1 +(1-\lambda)x_2 \in \setXinv$ and so $\lambda u_1 +(1-\lambda)u_2 \in \setU(\lambda x_1 +(1-\lambda)x_2)$. Therefore  $ \lambda (x_1,u_1) + (1-\lambda)(x_2,u_2) \in Gr(\setU)$, for any $\lambda\in[0,1]$, and so $Gr(\setU)$ is convex.
\end{proof}

\begin{lem}[Upper hemicontinuity of $\setU(x)$]\label{lem:uppercontin}
	Consider DICSys \eqref{ec:sysDisc} and a compact and convex CIS $\setXinv \subset \setX$.
	Then correspondence $\Ucal(x) := \{ u \in \setU : A^dx + B^du \in \setXinv\}$ is upper hemicontinuous.
\end{lem}

\begin{proof}
	The proof follows directly from Lemmas \ref{lem:uppcont}, \ref{lem:imagclosure} and \ref{lem:graphclosure}.
\end{proof}

In order to prove the lower hemicontinuity, let us first introduce some preliminary definition and result.

\begin{defi}[Radially convex set]
	A set $\setX$ is radially convex at a point $x\in\setX$ if there exists $\epsilon>0$ such that for every $z\in \setX\cap B(x,\epsilon)$, with $z\ne x$, there exists $w\in\setX\cap\partial B(x,\epsilon)$ and $\lambda\in (0,1)$ for which $z=\lambda x + (1-\lambda)w$.
\end{defi}

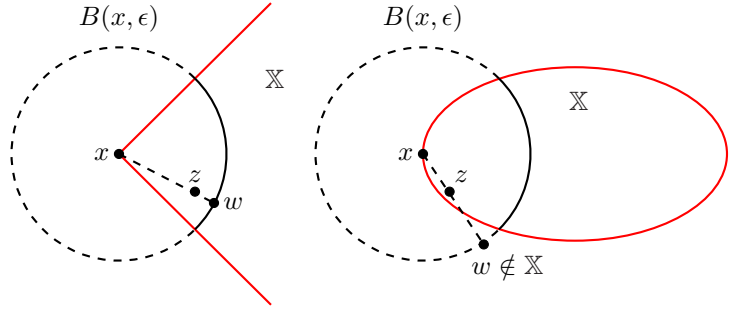
\begin{figure}[h!tb]
	\begin{center}
		\begin{tikzpicture}
		\draw[thick,color=red] (0,0)--(2,2);
		\draw[thick,color=red] (0,0)--(2,-2);
		\draw[thick] (1,-1) arc (-45:45:1.4142cm);
		\draw[thick,dashed] (1,1) arc (45:315:1.4142cm);
		\node[left] at (0,0) {\scalebox{1}{$x$}};
		\draw [fill=black](0,0) circle (0.06cm);
		\draw [fill=black](1,-0.5) circle (0.06cm);
		\node[above] at (1,-0.5) {\scalebox{1}{$z$}};
		\node[right] at (1.25,-0.65) {\scalebox{1}{$w$}};
		\draw [fill=black](1.25,-0.65) circle (0.06cm);
		\draw[thick,dashed] (0,0)--(1.25,-0.65) ;
		\node[right,color=red] at (1.8,1) {\scalebox{1}{$\setX$}};
		\node[above] at (0,1.5) {\scalebox{1}{$B(x,\epsilon)$}};
		\draw[thick,color=red] (6,0) ellipse (2cm and 1.15cm);
		\node[above] at (4,1.5) {\scalebox{1}{$B(x,\epsilon)$}};
		\draw[thick] (5,-1) arc (-45:45:1.4142cm);
		\draw[thick,dashed] (1+4,1) arc (45:315:1.4142cm);
		\node[left] at (0+4,0) {\scalebox{1}{$x$}};
		\draw [fill=black](0+4,0) circle (0.06cm);
		\draw [fill=black](4.35,-0.5) circle (0.06cm);
		\node[above] at (4.35,-0.5) {\scalebox{1}{\phantom{aa}$z$}};
		\node[right,below] at (4.8,-1.2) {\scalebox{1}{\phantom{aaaa}$w\notin\setX$}};
		\draw [fill=black](4.8,-1.2) circle (0.06cm);
		\draw[thick,dashed] (4,0)--(4.8,-1.2) ;
		\node[right,color=red] at (5.8,0.7) {\scalebox{1}{$\setX$}};
		\end{tikzpicture}
	\end{center}
	\caption{A radially convex set in $x$ (left) and one that it is not (right)}
	\label{fig:rad_conv_set}
\end{figure}

In words, a set $\setX$ is radially convex at a point $x\in\setX$ if there exists a ball such that every point of $\setX$ that is inside the ball, can be represented as a strict convex combination of the center of the ball and some point of $\setX$ that lies on the boundary of the ball. For example, in Figure~\ref{fig:rad_conv_set} we present two examples, one that is radially convex at a point $x$ and one that it is not.

Note that the radially convex set that was chosen in the figure as an example is a polytope. This is not casual, it is in indeed the general case in the context of compact convex sets. We state this result in the following lemma. 

\begin{lem}[Polytope characterization]\label{lem:poly_prop}
	Let $\setX$ be a compact convex set of $\mathbb R^n$, then the following assertion are equivalent:
	\begin{itemize}
		\item $\setX$ is polytope.
		\item $\setX$ is radially convex at every $x\in\setX$.
	\end{itemize}
\end{lem}

The proof of this result can be found in ~\cite[Lemma 1]{Mac06}. This result is the key stone to prove the lower hemicontinuity of $\setU(x)$ and brings out the necessity of working with polytopes.

\begin{lem}[Lower hemicontinuity of $\setU(x)$]\label{lem:lowercontin}
	Consider DICSys \eqref{ec:sysDisc} and a compact convex polytope CIS $\setXinv \subset \setX$.
	Then correspondence $\setU(x) := \{ u \in \setU : A^dx + B^du \in \setXinv\}$ is lower hemicontinuous.
\end{lem}

This result can also be found in~\cite[Theorem 2]{Mac06} in a more general context. For the sake of completeness and clarity we decided to include it here.

\begin{proof}
	Let $x\in\setXinv$ and $\{x_k\}\subset\setXinv$ a sequence such that $\{x_k\}\to x$, when $k\to\infty$. Let $u \in \setU(x)$. To show that the correspondence $\setU(x)$ is lower hemicontinuous we need to construct a sequence of controls $u_k\in\setU(x_k)$ such that $u_k\to u$, when $k\to\infty$.
	
	Since $\setXinv$ is a compact convex polytope by Lemma~\ref{lem:poly_prop} there exists $\epsilon>0$ such that $B(x,\epsilon)\cap\setXinv$ is radially convex at $x$. Since $\{x_k\}\to x$, there exists $K>0$ such that $x_k \in B(x,\epsilon)\cap\setXinv$ for every $k\ge K$. Then for every $k\ge K$ there exist $w_k\in\setXinv\cap\partial B(x,\epsilon)$ and $\lambda_k\in (0,1)$ for which $x_k=\lambda_k x + (1-\lambda_k)w_k$. Since $\{w_k\}$ is bounded ($w_k\in\setXinv$) and $x_k\to x$, we have that $\lambda_k \to 0$, when $k\to\infty$.
	
	For every $k\ge K$, fix a sequence of controls $v_k\in\setU(w_k)$. By Lemma~\ref{lem:graphclosure} $Gr(\setU)$ is convex, so $\lambda_k (x,u) + (1-\lambda_k)(w_k,v_k)\in Gr(\setU)$, i.e.  $(x_k,\lambda_ku+(1-\lambda_k)v_k) \in Gr(\setU)$. Hence $u_k:=\lambda_ku+(1-\lambda_k)v_k \in\setU(x_k)$. Since $\{v_k\}$ is bounded ($v_k\in\setU_{\setXinv}$) and $\lambda_k\to 0$, we have that $u_k \to u$ when $k\to\infty$. Therefore $\setU$ is lower hemicontinuous at~$x$.
\end{proof}

\section{Appendix C. Positive Polynomials on a Finite Interval}
For sake of self containment, we recall here some concepts introduced in sections 3.3, of \cite{Nesterov2000}.
The following sections discuss how the cone of cofficients of univariate polynomials can be represented as linear images of the cone of positive semi-definite matrices. This enables solving the corresponding optimization problems using semidefinite programming schemes.
Let $\mathcal{S} = \{u^{(1)}(x),\cdots,u^{(m)}(x)\}, x\in \Delta$, be an arbitrary system of linearly independent functions.
Definine the finite-dimensional functional subspace
\begin{equation*}
\mathcal{F}(\mathcal{S}) = \{q(x) = \sum_{k=1}^{m}q^{(k)}u^{(k)}(x), q = (q^{(1)},\cdots,q^{(m)}) \in \mathbb{R} \}.
\end{equation*}
The convex cone 
\begin{equation*}
K = \{p(x) = \sum_{i=1}^{N}q_i^2(x), q_i(x)\in \mathcal{F}(\mathcal{S}), i=1,\cdots,N\}
\end{equation*}
can be described by the squared functional system
\begin{equation*}
\mathcal{S}^2 = \{v_{ij}(x) = u^{(i)}(x)u^{(j)}(x), i,j=1,\cdots,m\}.
\end{equation*}
Define $v(x)$ the vector of components of a basis of $\mathcal{S}^2$.
Define the vector coefficients $\lambda_{ij} \in \mathbb{R}^n$ as:
\begin{equation*}
u^{(i)}(x) \cdot u^{(j)}(x) = \lambda_{ij}^T v(x), \forall x \in \Delta
\end{equation*}
A matrix valued operator can be defined as:
\begin{equation*}
(\Lambda(v))^{ij}u^{(i)}(x) \cdot u^{(j)}(x) = \lambda_{ij}^T v(x), \forall x \in \Delta
\end{equation*}
Note that 
$$u(x)u(x)^T \equiv \Lambda(v(x)), ~x \in \Delta$$
The adjoint linear operator $\Lambda^*(Y), Y \in \mathbb{R}^{m\times m}$ is defined as
$$ <Y, \Lambda(v)> = <\Lambda^*(Y),v> $$

\subsection{LMI representation of a cone}
\begin{theorem}
	\begin{enumerate}
		\item The function $p(x) = p^T v(x), p \in \mathbb{R}^n$, belongs to $K$ if and only if there exists a positive semidefinite $(m\times m)$-matrix $Y$ such that $p = \Lambda^*(Y)$:
	\begin{equation}
	K = \{p \in \mathbb{R}^n: p=\Lambda^*(Y), Y \succeq 0 \}.
	\end{equation}
	\item Any $p \in K$ can be represented as a sum of at most $m$ squares,
	\begin{equation*}
	p(x) = \sum_{i=1}^{k}q_i^2(x), q_i(x) \in \mathcal{F}(\mathcal{S}), i=1,\cdots,k\leq m.
	\end{equation*}
	\end{enumerate}
\end{theorem}

\subsection{Representation of a non-negative polynomial as LMI}
\label{sec:SoS}
Consider a fixed interval $[a,b] \subset \mathbb{R}$, the vector function $v(w) = (1,w,w^2,\cdots, w^n) \in \mathbb{R}^{n+1}, w \in [w,b]$.
Define the following convex cone:
\begin{equation*}
    K_\mathrm{a,b}={p \in \mathbb{R}^{n+1}:p^T v(w) \geq 0, \forall w \in [a,b]}.
\end{equation*}
This cone describes the set of coefficient vectors such that the polynomial $p^T v(w)$ is non-negative for all $w \in [a,b]$.
Now, let $n = 2 m$ and denote 
\begin{equation*}
u_1(w)=(1,w,\cdots,w^m) \in \mathbb{R}^{m+1},
\end{equation*}
\begin{equation*}
u_2(w)=(1,w,\cdots,w^m) \in \mathbb{R}^{m-1},
\end{equation*}
From Markov-Lukasz theorem, non negative polynomials can be represented as
\begin{equation}
p(w) = (q_1^T u_1(w))^2 + (w-a)(b-w)(q_2^T u_2(w))^2
\end{equation}
with some $q_1 \in  \mathbb{R}^{m+1},$ $q_2 \in  \mathbb{R}^{m}$. 
This configures a sum of weigthed squares.
\subsection{Matrix exponential}
Given a square $n \times n$ diagonalizable transition matrix $A = V \Lambda V^{-1}$ with eigensolutions $A x_k = \lambda_k x_k, (k=1,2,\cdots, n)$.
The solution to the system $\dot{x} = A x$ can be written as
\begin{equation*}
x(t) = \sum_{i=1}^{n}c_i e^{\lambda_i t}v_i
\end{equation*}
where $c = V^{-1} x(0)$. 
\subsection{Modal Coordinates}
Considering a dynamic system wint $n$ distinct eigenvalues, its free response can be decomposed in modal form
\begin{equation}
x(t) = e^{At} x(0) = e^{V \Lambda V^{-1} t} x(0),
\end{equation}
where the $i$-th column of $V = [v_1 \cdots v_n]$ is the $i$-th eigenvector and $\Lambda$ is a diagonal matrix with the eigenvalues as its entries.
Noticing that for any matrix B and nonsingular matrix P
\begin{equation}
e^{V \Lambda V^{-1} t} = V e^{\Lambda t} V^{-1},
\end{equation}
or, equivalently,
\begin{equation}
e^{A t} = V \begin{bmatrix}
e^{\lambda_1 t} & 0 & \cdots & 0 \\
0 & e^{\lambda_2 t} & \cdots & 0 \\
& & \cdots & 0\\
0 &\cdots &0 & e^{\lambda_n t}
\end{bmatrix} 
V^{-1},
\end{equation}
Moreover, this can be expressed as
\begin{equation}
e^{A t} = \sum_{i=1}^n V \delta_i^n
V^{-1} e^{\lambda_n t}
\end{equation}
or in a short form as
\begin{equation}
e^{A t} = \sum_{i=1}^n \phi_i e^{\lambda_n t}
\end{equation}
using $\phi_i = V \delta_i^n V^{-1}$ and $\delta_i^n$ is a $n \times n$ matrix filled with zeros except at the $i,i$-th element where it is 1.
\subsection{Spectral Representation of the Positive Definite Matrix}{\label{spectral-representation}}
Suppose positive semidefinite (symmetric) matrix $A \in \mathcal{R}^{n \times n}$, all eigenvalues $\lambda_i$ are distinct, for linear independent eigenvectors $a_i$, with normalized $a_i' a_j  = \delta_{ij}$. Decomposing $A = \sum_{i=1}^{n} a_i a_i' \lambda_i$ and applying the well-known Taylor expansion for $e^A$,
\begin{equation}
e^A = \sum_{i=1}^{n}  a_i a_i' e^\lambda_i
\end{equation}
\begin{equation}
\Phi(t) = e^{At} = \sum_{i=1}^{n}  a_i a_i' e^{\lambda_i t}
\end{equation}
denoting $\alpha^r = a_r a_r'$ the matrix corresponding to the $r$-th eigenvector and $\alpha^r_{ij}$ the $ij$-th singleton, each singleton in the transition matrix can be expressed as
\begin{equation}
\Phi_{ij}(t) = \sum_{r=1}^{n}  \alpha^r_{ij} e^{\lambda_r t}
\end{equation}

\bibliographystyle{plain} 
\bibliography{TrackImp,biblio_rdv,bib_impuls,Propios,approximation}

\end{document}